# CONVERGENCE RATE AND AVERAGING OF NONLINEAR TWO-TIME-SCALE STOCHASTIC APPROXIMATION ALGORITHMS

By Abdelkader Mokkadem and Mariane Pelletier

*University of Versailles–Saint-Quentin*

The first aim of this paper is to establish the weak convergence rate of nonlinear two-time-scale stochastic approximation algorithms. Its second aim is to introduce the averaging principle in the context of two-time-scale stochastic approximation algorithms. We first define the notion of asymptotic efficiency in this framework, then introduce the averaged two-time-scale stochastic approximation algorithm, and finally establish its weak convergence rate. We show, in particular, that both components of the averaged two-time-scale stochastic approximation algorithm simultaneously converge at the optimal rate $\sqrt{n}$.

**1. Introduction.** Let

$$f : \begin{cases} \mathbb{R}^d \times \mathbb{R}^{d'} \to \mathbb{R}^d \\ (\theta, \mu) \mapsto f(\theta, \mu) \end{cases} \quad \text{and} \quad g : \begin{cases} \mathbb{R}^d \times \mathbb{R}^{d'} \to \mathbb{R}^{d'} \\ (\theta, \mu) \mapsto g(\theta, \mu) \end{cases}$$

be two unknown functions, and let $(\theta^*, \mu^*)$ be the unique solution to the equations

$$f(\theta, \mu) = 0 \quad \text{and} \quad g(\theta, \mu) = 0.$$

Assume that error-contaminated observations of $f(\theta, \mu)$ and $g(\theta, \mu)$ are available at any level $(\theta, \mu)$. The two-time-scale stochastic approximation algorithm, which allows the recursive approximation of $(\theta^*, \mu^*)$, is defined as

$$\theta_{n+1} = \theta_n + \beta_n X_{n+1}, \tag{1}$$

$$\mu_{n+1} = \mu_n + \gamma_n Y_{n+1}, \tag{2}$$

where $X_{n+1}$ and $Y_{n+1}$ are error-contaminated observations of $f(\theta_n, \mu_n)$ and $g(\theta_n, \mu_n)$, respectively, and where the step sizes $(\beta_n)$ and $(\gamma_n)$ are two positive nonrandom sequences converging to zero with different rates.

Over the past few years, several such algorithms have been proposed for various applications (see [1, 3, 4, 12, 13]), and criteria ensuring the almost sure convergence of $(\theta_n, \mu_n)$ to $(\theta^*, \mu^*)$ have been established by Borkar [5], Konda and Borkar [12] and Konda and Tsitsiklis [13]. To our knowledge, the only existing result on the convergence rate of the two-time-scale stochastic approximation algorithm (1)–(2) is the one of Konda and Tsitsiklis [14]. In the case when the functions $f$ and $g$ are linear and when $\lim_{n\to\infty}\beta_n/\gamma_n = 0$, Konda and Tsitsiklis [14] establish that the fastest component $\theta_n$ satisfies the following central limit theorem (CLT):

$$(3) \qquad \sqrt{\beta_n^{-1}}(\theta_n - \theta^*) \xrightarrow{\mathcal{D}} \mathcal{N}(0, \Sigma_\theta),$$

where $\xrightarrow{\mathcal{D}}$ denotes the convergence in distribution, $\mathcal{N}$ the Gaussian-distribution, and where the asymptotic covariance matrix $\Sigma_\theta$ is defined in (8) below. Moreover, it can be conjectured from their analysis that the slowest component $\mu_n$ fulfills the CLT:

$$(4) \qquad \sqrt{\gamma_n^{-1}}(\mu_n - \mu^*) \xrightarrow{\mathcal{D}} \mathcal{N}(0, \Sigma_\mu)$$

[where the asymptotic covariance matrix $\Sigma_\mu$ is defined in (9) below]. The result (3) of [14] is thus very surprising. As a matter of fact, it shows that the slowest component $\mu_n$ [which, through $X_{n+1}$, is present in the recursive definition (1) of $\theta_n$] has no effect on the convergence rate of the fastest component $\theta_n$, except in the expression of the asymptotic covariance matrix $\Sigma_\theta$. It is then natural to wonder whether this phenomenon is specific to the case of the functions $f$ and $g$ being linear or not.

Our first aim in this paper is to study the weak joint convergence rate of $\theta_n$ and $\mu_n$ in the case where the functions $f$ and $g$ are nonlinear. We still consider the case $\lim_{n\to\infty}\beta_n/\gamma_n = 0$, and prove that

$$(5) \qquad \begin{pmatrix} \sqrt{\beta_n^{-1}}(\theta_n - \theta^*) \\ \sqrt{\gamma_n^{-1}}(\mu_n - \mu^*) \end{pmatrix} \xrightarrow{\mathcal{D}} \mathcal{N}\left(0, \begin{pmatrix} \Sigma_\theta & 0 \\ 0 & \Sigma_\mu \end{pmatrix}\right).$$

The CLT (5) extends, in particular, the result (3) of [14] to the case where the functions $f$ and $g$ are nonlinear. Let us underline that, as explained in [14], in the case $(\beta_n) \equiv (\gamma_n)$, the algorithm defined by (1)–(2) reduces to a single-time-scale stochastic approximation algorithm used for the search of the zero of the function $h : \mathbb{R}^{d+d'} \to \mathbb{R}^{d+d'}$ defined by $h(\theta, \mu) = (f(\theta, \mu), g(\theta, \mu))$. The convergence rate of such single-time-scale stochastic approximation algorithms has been widely studied (see, among many others, Nevel'son and Has'minskii [21], Kushner and Clark [15], Benveniste, Métivier and Priouret [2], Ljung, Pflug and Walk [10] and Duflo [9]), but the existing techniques do not apply when the step sizes $(\beta_n)$ and $(\gamma_n)$ have two different convergence rates, that is, in the context of two-time-scale stochastic approximation algorithms. Let us also point out that the two-time-scale iterations considered by Konda and Tsitsiklis [14] and in the present paper are totally different from those that arise in the study of the tracking ability of adaptative algorithms (see [2]) or in the joint approximation of the location and size of the maximum of a regression function (see [20]); the specific difficulty in the present context relies on the double dependency between both components $\theta_n$ and $\mu_n$ [$\theta_n$ defined by (1) depends, through $X_{n+1}$, on $\mu_n$ defined by (2), whereas $\mu_n$ defined by (2) depends, through $Y_{n+1}$, on $\theta_n$ defined by (1)]. Let us finally underline that the techniques we use to prove (5) (introduction of exponential martingales and recourse to successive almost sure upper bounds) radically differ from those employed by Konda and Tsitsiklis [14] to establish (3); let us also mention that the additional difficulty induced by the nonlinearity of the functions $f$ and $g$ will be enlightened in our proof of (5).

Now, let us note that, in view of (3) and (5), the recommended choice of the fastest step size $(\beta_n)$ is $(\beta_n) \equiv (\beta_0 n^{-1})$, since it is the choice which ensures that the fastest component $\theta_n$ converges with the optimal rate $\sqrt{n}$. However, this optimal choice induces conditions on the parameter $\beta_0$, which are difficult to handle because of depending on an unknown parameter. The problem due to the choice of the optimal step size $(\beta_n) \equiv (\beta_0 n^{-1})$ is now well known in the context of single-time-scale stochastic approximation algorithms, and the method widely employed in this framework to circumvent this problem is the use of the averaging principle independently introduced by Ruppert [26] and Polyak [24], and then widely discussed and extended (see, among many others, Yin [27], Delyon and Juditsky [6], Polyak and Juditsky [25], Kushner and Yang [16], Dippon and Renz [7, 8], Duflo [9], Kushner and Yin [17] and Pelletier [23]).

Our second aim in this paper is to introduce the averaging principle in the context of two-time-scale stochastic approximation algorithms. We first define the notion of asymptotic efficiency in this framework, then introduce the averaged two-time-scale stochastic approximation algorithm, and finally establish its weak convergence rate. We prove, in particular, that, by choosing the step sizes $(\beta_n)$ and $(\gamma_n)$ equal to $(\beta_n) \equiv (\beta_0 n^{-b})$ and $(\gamma_n) \equiv (\gamma_0 n^{-a})$ with $1/2 < a < b < 1$, and by defining the averaged two-time-scale algorithm by setting

$$\overline{\theta}_n = \frac{1}{n} \sum_{k=1}^{n} \theta_k \quad \text{and} \quad \overline{\mu}_n = \frac{1}{n} \sum_{k=1}^{n} \mu_k,$$

where $\theta_k$ and $\mu_k$ are defined in (1)–(2), we obtain an asymptotically efficient two-time-scale algorithm, which satisfies the CLT

$$\begin{pmatrix} \sqrt{n}(\overline{\theta}_n - \theta^*) \\ \sqrt{n}(\overline{\mu}_n - \mu^*) \end{pmatrix} \xrightarrow{\mathcal{D}} \mathcal{N}(0, C),$$

where the asymptotic covariance matrix $C$ is precisely defined (see Theorem 2). The striking aspect of this result is that averaging leads to a two-time-scale algorithm whose components $\overline{\theta}_n$ and $\overline{\mu}_n$ simultaneously converge with the optimal rate $\sqrt{n}$.

Our paper is now organized as follows. Section 2 is devoted to the study of the convergence rate of nonlinear two-time-scale stochastic approximation algorithms. We first precisely state our assumptions and main results; then, we give the outlines of the proof of our main results, postponing the technical parts until the Appendix. Section 3 is reserved for averaging. The notion of asymptotic efficiency of two-time-scale stochastic approximation algorithms is introduced in Section 3.1; the weak convergence rate of the averaged two-time-scale algorithm is stated and then proved in Sections 3.2 and 3.3, respectively.

## 2. Convergence rate of nonlinear two-time-scale stochastic approximation algorithms.

2.1. *Assumptions and notation.* For any square matrix $A$, we set

$$\Lambda^{(A)} = -\max\{\mathcal{R}e(\lambda), \lambda \in \mathrm{Sp}(A)\},$$

where $\mathrm{Sp}(A)$ denotes the spectrum of $A$. Moreover, $\|\cdot\|$ denotes the Euclidean vector norm in $\mathbb{R}^d$, $\mathbb{R}^{d'}$ and $\mathbb{R}^{d+d'}$ without distinction, and $\|\cdot\|$ the matrix norm induced by the Euclidean vector norm.

The assumptions we require are the following:

(A1) $\lim_{n\to\infty} \theta_n = \theta^*$ a.s. and $\lim_{n\to\infty} \mu_n = \mu^*$ a.s.

(A2)   (i) There exists a neighborhood $\mathcal{U}$ of $(\theta^*, \mu^*)$ such that, for all $(\theta, \mu) \in \mathcal{U}$,

$$\begin{pmatrix} f(\theta, \mu) \\ g(\theta, \mu) \end{pmatrix} = \begin{pmatrix} Q_{11} & Q_{12} \\ Q_{21} & Q_{22} \end{pmatrix} \begin{pmatrix} \theta - \theta^* \\ \mu - \mu^* \end{pmatrix} + O\left( \left\| \begin{matrix} \theta - \theta^* \\ \mu - \mu^* \end{matrix} \right\|^2 \right).$$

(ii) Set

(6) $$H = Q_{11} - Q_{12} Q_{22}^{-1} Q_{21}.$$

We have $\Lambda^{(H)} > 0$ and $\Lambda^{(Q_{22})} > 0$.

(A3)   (i) $(\beta_n) \equiv (\beta_0 n^{-b})$ and $(\gamma_n) \equiv (\gamma_0 n^{-a})$ with $\beta_0 > 0$, $\gamma_0 > 0$ and $\frac{1}{2} < a < b \leq 1$.

(ii) If $b = 1$, then $\beta_0 > 1/[2\Lambda^{(H)}]$.

(A4) The error-contaminated observations can be written as

$$X_{n+1} = f(\theta_n, \mu_n) + \psi_n^{(\theta)} + V_{n+1},$$

$$Y_{n+1} = g(\theta_n, \mu_n) + \psi_n^{(\mu)} + W_{n+1},$$

and denoting by $\mathcal{F}_n$ the $\sigma$-field spanned by $\{V_i, W_j, \theta_k, \mu_l, \psi_{k'}^{(\theta)}, \psi_{l'}^{(\mu)}, 0 \leq i, j, k, l, k', l' \leq n\}$, we have the following:

(i) $\mathbb{E}(V_{n+1} | \mathcal{F}_n) = 0$ and $\mathbb{E}(W_{n+1} | \mathcal{F}_n) = 0$ a.s.

(ii) There exists a positive matrix $\Gamma$ such that

$$\lim_{n \to \infty} \mathbb{E}\left( \begin{pmatrix} V_{n+1} \\ W_{n+1} \end{pmatrix} (V_{n+1}^T \ W_{n+1}^T) \Big| \mathcal{F}_n \right) = \Gamma = \begin{pmatrix} \Gamma_{11} & \Gamma_{12} \\ \Gamma_{21} & \Gamma_{22} \end{pmatrix} \qquad \text{a.s.}$$

(iii) There exists $m > 2/a$ such that $\sup_n \mathbb{E}(\|V_{n+1}\|^m | \mathcal{F}_n) < \infty$ and $\sup_n \mathbb{E}(\|W_{n+1}\|^m | \mathcal{F}_n) < \infty$ a.s.

(iv)

$$\psi_n^{(\theta)} = r_n^{(\theta)} + O(\|\theta_n - \theta^*\|^2 + \|\mu_n - \mu^*\|^2),$$

$$\psi_n^{(\mu)} = r_n^{(\mu)} + O(\|\theta_n - \theta^*\|^2 + \|\mu_n - \mu^*\|^2),$$

with $\|r_n^{(\theta)}\| + \|r_n^{(\mu)}\| = o(\sqrt{\beta_n})$ a.s.

Let us specify that the matrices $Q_{11}$ and $\Gamma_{11}$ (resp. $Q_{22}$ and $\Gamma_{22}$) in (A2)(i) and (A4)(ii) are $d \times d$ (resp. $d' \times d'$) matrices; the matrices $Q_{12}$, $Q_{21}$, $\Gamma_{12}$ and $\Gamma_{21}$ are of appropriate dimension. Set

$$(7) \quad \begin{aligned} \Gamma_\theta &= \lim_{n \to \infty} \mathbb{E}([V_{n+1} - Q_{12}Q_{22}^{-1}W_{n+1}][V_{n+1} - Q_{12}Q_{22}^{-1}W_{n+1}]^T | \mathcal{F}_n) \\ &= \Gamma_{11} + Q_{12}Q_{22}^{-1}\Gamma_{22}[Q_{22}^{-1}]^T Q_{12}^T - \Gamma_{12}[Q_{22}^{-1}]^T Q_{12}^T - Q_{12}Q_{22}^{-1}\Gamma_{21}. \end{aligned}$$

We can now give the explicit definition of the asymptotic covariance matrices $\Sigma_\theta$ and $\Sigma_\mu$, which stand in (3), (4) and (5):

$$(8) \qquad \Sigma_\theta = \int_0^\infty \exp\left[ \left( H + \frac{\mathbb{1}_{b=1}}{2\beta_0} I \right) t \right] \Gamma_\theta \exp\left[ \left( H^T + \frac{\mathbb{1}_{b=1}}{2\beta_0} I \right) t \right] dt,$$

$$(9) \qquad \Sigma_\mu = \int_0^\infty \exp[Q_{22}t] \Gamma_{22} \exp[Q_{22}t] \, dt.$$

Let us mention that the matrices $\Sigma_\theta$ and $\Sigma_\mu$ are the solutions of the Lyapounov equations

$$\left[ H + \frac{\mathbb{1}_{b=1}}{2\beta_0} I \right] \Sigma_\theta + \Sigma_\theta \left[ H^T + \frac{\mathbb{1}_{b=1}}{2\beta_0} I \right] = -\Gamma_\theta$$

and

$$Q_{22}\Sigma_\mu + \Sigma_\mu Q_{22}^T = -\Gamma_{22},$$

respectively (see Lemma 3.I.3 in (year?)).

*Comments on the assumptions.*

1. We refer to [5, 12, 13] for quite general conditions that ensure the consistency assumption (A1). Let us underline that, in the case where $f$ and $g$ are linear, $\psi_n^{(\theta)} = 0$ and $\psi_n^{(\mu)} = 0$, assumption (A1) is useless; as a matter of fact, as noted by Konda and Tsitsiklis [14], assumptions (A2)–(A4) imply (A1) in this particular case. Let us also mention that a particular example of two-time-scale stochastic approximation algorithm is the well known Polyak–Ruppert averaging; in this framework, (1)–(2) reduces to

$$\theta_{n+1} = \theta_n + \frac{1}{n}(\mu_n - \theta_n),$$

$$\mu_{n+1} = \mu_n + \gamma_n Y_{n+1},$$

where $Y_{n+1}$ is an error-contaminated observation at $\mu_n$ of an unknown function $h$, and $\lim_{n \to \infty} n\gamma_n = \infty$; (A1) then comes down to the assumption $\lim_{n \to \infty} \mu_n = \mu^*$ [where $h(\mu^*) = 0$], and conditions which ensure this lattest assumption can be found, among many others, in [9, 15, 18].

2. Assumptions (A2)(ii) and (A3)(ii) ensure that the matrices $\Sigma_\theta$ and $\Sigma_\mu$ are well defined. As a matter of fact, the conditions in (A2)(ii) mean that the matrices $H$ and $Q_{22}$ are attractive (or Hurwitz) and, in the case $b = 1$, it follows from the condition in (A3)(ii) that the matrix $[H + \frac{1}{2\beta_0}I]$ is attractive.

3. To establish the convergence rate of the two-time-scale stochastic approximation algorithm (1)–(2), Konda and Tsitsiklis [14] assume that the functions $f$ and $g$ are linear, that is, that

$$\begin{pmatrix} f(\theta,\mu) \\ g(\theta,\mu) \end{pmatrix} = \begin{pmatrix} Q_{11} & Q_{12} \\ Q_{21} & Q_{22} \end{pmatrix} \begin{pmatrix} \theta - \theta^* \\ \mu - \mu^* \end{pmatrix}.$$

Moreover, their framework corresponds to the case (A4) is fulfilled with $\psi_n^{(\theta)} = 0$, $\psi_n^{(\mu)} = 0$, and $(V_n, W_n)$ are independent random vectors with zero mean and common covariance $\Gamma$. On the other hand, their conditions on the step sizes $(\beta_n)$ and $(\gamma_n)$ are more general than ours.

2.2. *Main results.* Our main result in this section is the following theorem.

THEOREM 1 [Joint weak convergence rate of $(\theta_n)$ and $(\mu_n)$]. *Let $(\theta_n, \mu_n)$ be defined by the recursive equations* (1)–(2). *Under assumptions* (A1)–(A4), *we have*

$$\begin{pmatrix} \sqrt{\beta_n^{-1}}(\theta_n - \theta^*) \\ \sqrt{\gamma_n^{-1}}(\mu_n - \mu^*) \end{pmatrix} \xrightarrow{\mathcal{D}} \mathcal{N}\left(0, \begin{pmatrix} \Sigma_\theta & 0 \\ 0 & \Sigma_\mu \end{pmatrix}\right),$$

*where $\Sigma_\theta$ and $\Sigma_\mu$ are defined in* (8) *and* (9), *respectively.*

The following proposition, which is of independent interest, will be a key tool for the study of the weak convergence rate of the averaged two-time-scale stochastic approximation algorithm.

PROPOSITION 1 [Strong convergence rate of $(\theta_n)$ and $(\mu_n)$]. *Let* $(\theta_n, \mu_n)$ *be defined by the recursive equations* (1)–(2). *Under assumptions* (A1)–(A4), *we have*

$$\|\theta_n - \theta^*\| = O\left(\sqrt{\beta_n \log\left[\sum_{k=1}^{n} \beta_k\right]}\right) \qquad a.s.$$

*and*

$$\|\mu_n - \mu^*\| = O\left(\sqrt{\gamma_n \log\left[\sum_{k=1}^{n} \gamma_k\right]}\right) \qquad a.s.$$

2.3. *Proof of Theorem* 1 *and Proposition* 1. Throughout the proof of Theorem 1 and Proposition 1, we assume, without loss of generality, that $\theta^* = 0$ and $\mu^* = 0$. In view of assumptions (A1), (A2) and (A4), we can write

$$(10) \qquad \theta_{n+1} = \theta_n + \beta_n(Q_{11}\theta_n + Q_{12}\mu_n + \rho_n^{(\theta)} + r_n^{(\theta)} + V_{n+1}),$$

$$(11) \qquad \mu_{n+1} = \mu_n + \gamma_n(Q_{21}\theta_n + Q_{22}\mu_n + \rho_n^{(\mu)} + r_n^{(\mu)} + W_{n+1}),$$

where

$$(12) \qquad \|\rho_n^{(\theta)}\| = O(\|\theta_n\|^2 + \|\mu_n\|^2) \quad \text{and} \quad \|\rho_n^{(\mu)}\| = O(\|\theta_n\|^2 + \|\mu_n\|^2).$$

Note that (11) gives

$$\mu_n = Q_{22}^{-1}\gamma_n^{-1}[\mu_{n+1} - \mu_n] - Q_{22}^{-1}(Q_{21}\theta_n + \rho_n^{(\mu)} + r_n^{(\mu)} + W_{n+1}),$$

and thus, in view of (10), it follows that

$$
\begin{aligned}
\theta_{n+1} &= \theta_n + \beta_n(Q_{11}\theta_n + Q_{12}Q_{22}^{-1}\gamma_n^{-1}[\mu_{n+1} - \mu_n] \\
&\qquad - Q_{12}Q_{22}^{-1}(Q_{21}\theta_n + \rho_n^{(\mu)} + r_n^{(\mu)} + W_{n+1}) \\
&\qquad + \rho_n^{(\theta)} + r_n^{(\theta)} + V_{n+1}) \\
&= \theta_n + \beta_n H\theta_n + \beta_n Q_{12}Q_{22}^{-1}\gamma_n^{-1}[\mu_{n+1} - \mu_n] \\
&\qquad + \beta_n(V_{n+1} - Q_{12}Q_{22}^{-1}W_{n+1}) \\
&\qquad + \beta_n([\rho_n^{(\theta)} + r_n^{(\theta)}] - Q_{12}Q_{22}^{-1}[\rho_n^{(\mu)} + r_n^{(\mu)}]),
\end{aligned}
$$

(13)

where $H$ is defined in (6). Now, set

$$u_n = \sum_{k=1}^{n} \beta_k,$$

$$(14) \qquad L_{n+1}^{(\theta)} = e^{u_n H} \sum_{k=1}^n e^{-u_k H} \beta_k (V_{k+1} - Q_{12} Q_{22}^{-1} W_{k+1}),$$

$$(15) \qquad R_{n+1}^{(\theta)} = e^{u_n H} \sum_{k=1}^n e^{-u_k H} \beta_k Q_{12} Q_{22}^{-1} \gamma_k^{-1} [\mu_{k+1} - \mu_k],$$

$$(16) \qquad \Delta_{n+1}^{(\theta)} = \theta_{n+1} - L_{n+1}^{(\theta)} - R_{n+1}^{(\theta)}$$

and

$$s_n = \sum_{k=1}^n \gamma_k,$$

$$(17) \qquad L_{n+1}^{(\mu)} = e^{s_n Q_{22}} \sum_{k=1}^n e^{-s_k Q_{22}} \gamma_k W_{k+1},$$

$$(18) \qquad R_{n+1}^{(\mu)} = e^{s_n Q_{22}} \sum_{k=1}^n e^{-s_k Q_{22}} \gamma_k Q_{21} [L_k^{(\theta)} + R_k^{(\theta)}],$$

$$(19) \qquad \Delta_{n+1}^{(\mu)} = \mu_{n+1} - L_{n+1}^{(\mu)} - R_{n+1}^{(\mu)}.$$

The main idea to establish Theorem 1 and Proposition 1 is to prove that the sequences $(R_n^{(\theta)})$ and $(\Delta_n^{(\theta)})$ are negligible in front of $(L_n^{(\theta)})$ on the one hand, and that the sequences $(R_n^{(\mu)})$ and $(\Delta_n^{(\mu)})$ are negligible in front of $(L_n^{(\mu)})$ on the other hand; the convergence rates of $(\theta_n)$ and $(\mu_n)$ are then given by the ones of $(L_n^{(\theta)})$ and $(L_n^{(\mu)})$, respectively. Let us note that, even though the sequence $(\mu_n)$ goes to zero a.s. slower than the sequence $(\theta_n)$ does, we shall prove that the term $(R_n^{(\theta)})$ goes to zero a.s. faster than the sequence $(\theta_n)$ does. This is due to an averaging effect, the sequence $(R_n^{(\theta)})$ bringing in a weighted sum of the differences $\mu_{k+1} - \mu_k$. In the sequel we shall come back on this effect several times.

Applying Lyapounov's theorem, we obtain the following lemma (see Section A.2 for the technical details).

LEMMA 1 [Joint weak convergence rate of $(L_n^{(\theta)})$ and $(L_n^{(\mu)})$].  *We have*

$$\begin{pmatrix} \sqrt{\beta_n^{-1}} L_n^{(\theta)} \\ \sqrt{\gamma_n^{-1}} L_n^{(\mu)} \end{pmatrix} \xrightarrow{\mathcal{D}} \mathcal{N} \left( 0, \begin{pmatrix} \Sigma_\theta & 0 \\ 0 & \Sigma_\mu \end{pmatrix} \right).$$

Moreover, the following lemma is proved in [22].

LEMMA 2 [Strong convergence rate of $(L_n^{(\theta)})$ and $(L_n^{(\mu)})$].  *We have*

$$\|L_n^{(\theta)}\| = O(\sqrt{\beta_n \log u_n}) \qquad a.s.$$

*and*

$$\|L_n^{(\mu)}\| = O(\sqrt{\gamma_n \log s_n}) \qquad a.s.$$

Theorem 1 (resp. Proposition 1) thus follows from the combination of Lemma 1 (resp. of Lemma 2), and of the following two lemmas (which imply, in particular, that the sequences $(\beta_n^{-1/2}[R_n^{(\theta)} + \Delta_n^{(\theta)}])$ and $(\gamma_n^{-1/2}[R_n^{(\mu)} + \Delta_n^{(\mu)}])$ go to zero a.s.):

LEMMA 3 [Strong convergence rate of $(R_n^{(\theta)})$ and $(R_n^{(\mu)})$].

1. *There exists $s > b/2$ such that $\|R_n^{(\theta)}\| = O(n^{-s})$ a.s.*
2. $\|R_n^{(\mu)}\| = O(\sqrt{\beta_n \log u_n})$ *a.s.*

LEMMA 4 [Strong convergence rate of $(\Delta_n^{(\theta)})$ and $(\Delta_n^{(\mu)})$]. *We have*

$$\|\Delta_n^{(\theta)}\| = o(\sqrt{\beta_n}) \qquad a.s.,$$
$$\|\Delta_n^{(\mu)}\| = o(\sqrt{\beta_n}) \qquad a.s.$$

The key point in the proof of Theorem 1 and Proposition 1 is thus the proof of Lemmas 3 and 4. The rest of Section 2 is devoted to this proof (we shall refer to the Appendix for the technical details). Let us first give the strategy to prove these lemmas.

We note that, to obtain an upper bound of $(R_n^{(\mu)})$, we need to have an upper bound of $(R_n^{(\theta)})$, which requires to have an upper bound of $(\mu_n)$. The main idea to prove Lemma 3 is thus to proceed by successive upper bounds. In a first step, we shall start with the only upper bound of $(\mu_n)$ available to us, that is, in view of assumption (A1), with $\|\mu_n\| = o(1)$. This will enable us to establish a first upper bound of $(R_n^{(\theta)})$ and then of $(R_n^{(\mu)})$. With these preliminary upper bounds, we shall be able to prove preliminary upper bounds for $(\Delta_n^{(\theta)})$ and $(\Delta_n^{(\mu)})$. Using (19) and applying Lemma 2, we shall then slightly improve the first upper bound of $(\mu_n)$; starting with this second upper bound of $(\mu_n)$, we shall then repeat the procedure previously described to find a third upper bound of $(\mu_n)$, which slightly improves the second one, and we shall carry on these successive upper bounds until we obtain the adequate upper bounds of $(\mu_n)$, $(R_n^{(\theta)})$, $(R_n^{(\mu)})$, $(\Delta_n^{(\theta)})$ and $(\Delta_n^{(\mu)})$.

Let us mention that the step, which consists in deducing upper bounds of $(\Delta_n^{(\theta)})$ and $(\Delta_n^{(\mu)})$ from upper bounds of $(\mu_n)$, $(L_n^{(\theta)})$, $(L_n^{(\mu)})$, $(R_n^{(\theta)})$ and $(R_n^{(\mu)})$, is quite straightforward in the case when the functions $f$ and $g$ are linear, $\psi_n^{(\theta)} = 0$ and $\psi_n^{(\mu)} = 0$ (see Remark 4 below); however, in the case where the functions $f$ and $g$ are nonlinear, this step too requires to compute successive upper bounds [we shall first show that $\|\Delta_n^{(\mu)}\| = o(1)$, and then shall recursively improve the upper bound of $(\Delta_n^{(\mu)})$ until we find the adequate upper bound of $(\Delta_n^{(\mu)})$].

Our proof of Lemmas 3 and 4 is now organized as follows. We first define Conditions (C) and (C′) [that are expressed with respect to the step sizes $(\beta_n)$ and $(\gamma_n)$ resp.] for a nonrandom sequence, conditions which will be used throughout the proof. Then, in Section 2.3.1, we show how the knowledge of an upper bound of $(\mu_n)$ and of $(\Delta_n^{(\mu)})$ enables to establish upper bounds of $(R_n^{(\theta)})$, $(R_n^{(\mu)})$, $(\Delta_n^{(\theta)})$, and to improve the upper bound of $(\Delta_n^{(\mu)})$. Section 2.3.2 is devoted to the body of the proof of Lemmas 3 and 4.

DEFINITION 1 [Condition (C)].    Let $(w_n)$ be a sequence of real numbers. We say that $(w_n)$ satisfies Condition (C) if $(w_n)$ is positive and bounded and if:

- in the case $b = 1$, there exist $\omega \geq 0$ and a nondecreasing slowly varying function $\mathcal{L}$ such that $w_n = n^{-\omega}\mathcal{L}(n)$;
- in the case $b < 1$,
$$\frac{w_n}{w_{n+1}} = 1 + o(\beta_n).$$

DEFINITION 2 [Condition (C′)].    Let $(w_n)$ be a sequence of real numbers. We say that $(w_n)$ satisfies Condition (C′) if $(w_n)$ is positive and bounded and if
$$\frac{w_n}{w_{n+1}} = 1 + o(\gamma_n).$$

REMARK 1.    If $b = 1$ and if $(w_n)$ satisfies Condition (C) with $\omega = 0$, then the function $\mathcal{L}$ is necessary bounded.

REMARK 2.    In the case $b < 1$, if $(w_n)$ satisfies Condition (C), then $(w_n)$ satisfies Condition (C′).

2.3.1. *Intermediate upper bounds.*    We can now state the following lemma, which gives an upper bound of $(R_n^{(\theta)})$ and $(R_n^{(\mu)})$ under the assumption $\|\mu_n\| = O(w_n)$, where $(w_n)$ is a nonrandom sequence satisfying Conditions (C) and (C′). The proof of this lemma only requires classical computations, and is thus postponed until the Appendix (see Section A.3).

LEMMA 5 [Intermediate upper bound of $(R_n^{(\theta)})$ and $(R_n^{(\mu)})$].    *Assume that there exists a nonrandom sequence $(w_n)$ satisfying Conditions (C) and (C′),* *and such that* $\|\mu_n\| = O(w_n)$ *a.s. For all* $s \in ]1/2, \beta_0 \Lambda^{(H)}[$, *we have*

$$\|R_n^{(\theta)}\| = O(\beta_n \gamma_n^{-1} w_n + n^{-s}) \qquad a.s.,$$

$$\|R_n^{(\mu)}\| = O(\beta_n \gamma_n^{-1} w_n + \sqrt{\beta_n \log u_n}) \qquad a.s.$$

REMARK 3. The term $R_n^{(\mu)}$ can be seen as a (matricial) weighted average of the terms $L_k^{(\theta)} + R_k^{(\theta)}$; the second upper bound in Lemma 5 is established by proving that the same upper bound holds for the sequence $(L_n^{(\theta)} + R_n^{(\theta)})$ and for its average $(R_n^{(\mu)})$, which seems quite natural. On the other hand, the term $R_n^{(\theta)}$ can be seen as a (matricial) weighted average of the terms $\gamma_k^{-1}[\mu_{k+1} - \mu_k]$; the striking aspect of the first upper bound in Lemma 5 is that, although $\mu_n$ is bounded by $w_n$, although $\gamma_n^{-1} \to \infty$, the average $R_n^{(\theta)}$ can be bounded by $\beta_n \gamma_n^{-1} w_n$ (which is smaller than $w_n$ since $\beta_n \gamma_n^{-1} \to 0$). This averaging effect is similar to the one which appears in the study of the averaged single-time-scale stochastic approximation algorithm introduced by Ruppert [26] and Polyak [24].

We now state a lemma, which gives an upper bound of $(\Delta_n^{(\theta)})$ and $(\Delta_n^{(\mu)})$ under the assumption $\|\mu_n\| = O(w_n)$ and $\|\Delta_n^{(\mu)}\| = O(\delta_n^{(\mu)})$, where $(w_n)$ and $(\delta_n^{(\mu)})$ are two nonrandom sequences satisfying Conditions (C) and (C′).

LEMMA 6 [Intermediate upper bound of $(\Delta_n^{(\theta)})$ and $(\Delta_n^{(\mu)})$]. *Assume that there exist two nonrandom sequences* $(w_n)$ *and* $(\delta_n^{(\mu)})$ *satisfying Conditions* (C) *and* (C′), *and such that* $\|\mu_n\| = O(w_n)$ *a.s. and* $\|\Delta_n^{(\mu)}\| = O(\delta_n^{(\mu)})$ *a.s. We have*

$$\|\Delta_n^{(\theta)}\| = O(\beta_n^2 \gamma_n^{-2} w_n^2 + \beta_n \gamma_n^{-1} \delta_n^{(\mu)}) + o(\sqrt{\beta_n}) \qquad a.s.,$$

$$\|\Delta_n^{(\mu)}\| = O(\beta_n^2 \gamma_n^{-2} w_n^2 + \beta_n \gamma_n^{-1} \delta_n^{(\mu)}) + o(\sqrt{\beta_n}) \qquad a.s.$$

We now give the outlines of the proof of Lemma 6, and refer to the Appendix for the technical computations.

*Outlines of the proof of Lemma* 6. We first note that $\Delta_n^{(\theta)}$ and $\Delta_n^{(\mu)}$ satisfy the following recursive expressions (see Section A.4.1 for the algebra leading to these equations):

$$(20) \quad \begin{aligned} \Delta_{n+1}^{(\theta)} &= (I + \beta_n H)\Delta_n^{(\theta)} + O(\beta_n^2)[L_n^{(\theta)} + R_n^{(\theta)}] \\ &\quad + \beta_n([\rho_n^{(\theta)} + r_n^{(\theta)}] - Q_{12}Q_{22}^{-1}[\rho_n^{(\mu)} + r_n^{(\mu)}]), \end{aligned}$$

$$(21) \quad \begin{aligned} \Delta_{n+1}^{(\mu)} &= (I + \gamma_n Q_{22})\Delta_n^{(\mu)} + O(\gamma_n^2)[L_n^{(\mu)} + R_n^{(\mu)}] \\ &\quad + \gamma_n[\rho_n^{(\mu)} + r_n^{(\mu)} + Q_{21}\Delta_n^{(\theta)}]. \end{aligned}$$

Now, set $T$ and $M$ such that $\frac{\mathbb{1}_{b=1}}{2\beta_0} < T < \Lambda^{(H)}$ and $0 < M < \Lambda^{(Q_{22})}$ respectively. In view of Proposition 3.1.2 in [9], there exist two matrix norms $\|\!\|\cdot\|\!\|_T$ and $\|\!\|\cdot\|\!\|_M$, and there exists $a \in ]0, \inf\{1/T, 1/M\}[$ such that, for all $\gamma \leq a$, $\|\!\|I + \gamma H\|\!\|_T \leq 1 - \gamma T$ and $\|\!\|I + \gamma Q_{22}\|\!\|_M \leq 1 - \gamma M$. For $x$ in $\mathbb{R}^d$ (resp. in $\mathbb{R}^{d'}$), define $M^d(x) = [xx\cdots x]$ (resp. $M^{d'}(x) = [xx\cdots x]$) the $d \times d$ (resp. $d' \times d'$) matrix all of whose columns are $x$. The function $\|\cdot\|_T$ (resp. $\|\cdot\|_M$) defined on $\mathbb{R}^d$ (resp. on $\mathbb{R}^{d'}$) by $\|x\|_T = \|\!\|M^d(x)\|\!\|_T$ (resp. by $\|x\|_M = \|\!\|M^{d'}(x)\|\!\|_M$) is then a vector norm compatible with the matrix norm $\|\!\|\cdot\|\!\|_T$ (resp. with $\|\!\|\cdot\|\!\|_M$) (see [11], page 297). For $n$ large enough, we thus have

$$
\begin{aligned}
(22) \quad \|\Delta_{n+1}^{(\theta)}\|_T &\leq (1 - \beta_n T)\|\Delta_n^{(\theta)}\|_T + \beta_n[O(\beta_n\|L_n^{(\theta)}\|_T + \beta_n\|R_n^{(\theta)}\|_T)] \\
&\quad + \beta_n[O(\|\rho_n^{(\theta)}\|_T + \|r_n^{(\theta)}\|_T + \|Q_{12}Q_{22}^{-1}\rho_n^{(\mu)}\|_T \\
&\quad + \|Q_{12}Q_{22}^{-1}r_n^{(\mu)}\|_T)]
\end{aligned}
$$

and

$$
\begin{aligned}
(23) \quad \|\Delta_{n+1}^{(\mu)}\|_M &\leq (1 - \gamma_n M)\|\Delta_n^{(\mu)}\|_M + \gamma_n[O(\gamma_n\|L_n^{(\mu)}\|_M + \gamma_n\|R_n^{(\mu)}\|_M)] \\
&\quad + \gamma_n[O(\|\rho_n^{(\mu)}\|_M + \|r_n^{(\mu)}\|_M + \|Q_{21}\Delta_n^{(\theta)}\|_M)].
\end{aligned}
$$

REMARK 4. In the case where the functions $f$ and $g$ are linear and when $\psi_n^{(\theta)} = 0$ and $\psi_n^{(\mu)} = 0$, the terms $\rho_n^{(\theta)}$ and $\rho_n^{(\mu)}$ equal zero; replacing in (22) $\|L_n^{(\theta)}\|_T$ and $\|R_n^{(\theta)}\|_T$ by their upper bounds given in Lemmas 2 and 5 enables to get an upper bound $\delta_n^{(\theta)}$ of $\|\Delta_n^{(\theta)}\|_T$. Then, replacing in (23) $\|L_n^{(\mu)}\|_M$ and $\|R_n^{(\mu)}\|_M$ by their upper bounds given in Lemmas 2 and 5, and $\|Q_{21}\Delta_n^{(\theta)}\|_M$ by its upper bound $\delta_n^{(\theta)}$, enables to obtain an upper bound of $\|\Delta_n^{(\mu)}\|_M$. Thus, in this particular framework, the proof of intermediate upper bounds of $(\Delta_n^{(\theta)})$ and $(\Delta_n^{(\mu)})$ is quite straightforward. Moreover, the upper bounds of $(\Delta_n^{(\theta)})$ and $(\Delta_n^{(\mu)})$ obtained in this case are better than those stated in Lemma 6 [compare (22) with (28) below, and (23) with (27) below]; in particular, the knowledge of a preliminary upper bound $\delta_n^{(\mu)}$ of the sequence $(\Delta_n^{(\mu)})$ is not necessary.

By using the equivalence property of the finite-dimensional vector norms, we note that, in view of (12), (16) and (19), we have

$$
\begin{aligned}
\|\rho_n^{(\theta)}\|_T &+ \|Q_{12}Q_{22}^{-1}\rho_n^{(\mu)}\|_T \\
&= O(\|\rho_n^{(\theta)}\| + \|\rho_n^{(\mu)}\|)
\end{aligned}
$$

$$= O(\|\theta_n\|^2 + \|\mu_n\|^2)$$
$$= O(\|L_n^{(\theta)}\|^2 + \|R_n^{(\theta)}\|^2 + \|\Delta_n^{(\theta)}\|^2 + \|L_n^{(\mu)}\|^2 + \|R_n^{(\mu)}\|^2 + \|\Delta_n^{(\mu)}\|^2)$$
$$= O(\|L_n^{(\theta)}\|^2 + \|R_n^{(\theta)}\|^2 + \|\Delta_n^{(\theta)}\|_T^2 + \|L_n^{(\mu)}\|^2 + \|R_n^{(\mu)}\|^2 + \|\Delta_n^{(\mu)}\|_M^2).$$

It thus follows from (22) that there exists $C_1 > 0$ such that, for $n$ large enough,

$$
\begin{aligned}
(24) \quad \|\Delta_{n+1}^{(\theta)}\|_T \leq{} & (1 - \beta_n T)\|\Delta_n^{(\theta)}\|_T \\
& + \beta_n[O(\beta_n\|L_n^{(\theta)}\| + \beta_n\|R_n^{(\theta)}\| + \|r_n^{(\theta)}\| + \|r_n^{(\mu)}\|)] \\
& + \beta_n C_1(\|L_n^{(\theta)}\|^2 + \|R_n^{(\theta)}\|^2 + \|\Delta_n^{(\theta)}\|_T^2 \\
& \quad + \|L_n^{(\mu)}\|^2 + \|R_n^{(\mu)}\|^2 + \|\Delta_n^{(\mu)}\|_M^2).
\end{aligned}
$$

Similarly, we can deduce from (23) the existence of $C_2 > 0$ such that, for $n$ large enough,

$$
\begin{aligned}
(25) \quad \|\Delta_{n+1}^{(\mu)}\|_M \leq{} & (1 - \gamma_n M)\|\Delta_n^{(\mu)}\|_M \\
& + \gamma_n[O(\gamma_n\|L_n^{(\theta)}\| + \gamma_n\|R_n^{(\mu)}\| + \|r_n^{(\mu)}\|)] \\
& + \gamma_n C_2(\|L_n^{(\theta)}\|^2 + \|R_n^{(\theta)}\|^2 + \|\Delta_n^{(\theta)}\|_T^2 + \|L_n^{(\mu)}\|^2 + \|R_n^{(\mu)}\|^2 \\
& \quad + \|\Delta_n^{(\mu)}\|_M^2 + \|\Delta_n^{(\theta)}\|_T).
\end{aligned}
$$

Now, let us note that, in view of assumption (A1), we have $\lim_{n\to\infty}\theta_n = 0$ and $\lim_{n\to\infty}\mu_n = 0$ a.s. Since $\lim_{n\to\infty}\beta_n\gamma_n^{-1} = 0$, Lemma 5 [applied with the sequence $(w_n) \equiv 1$] implies that $\lim_{n\to\infty}R_n^{(\theta)} = 0$ and $\lim_{n\to\infty}R_n^{(\mu)} = 0$ a.s. Noting that Lemma 2 ensures that $\lim_{n\to\infty}L_n^{(\theta)} = 0$ and $\lim_{n\to\infty}L_n^{(\mu)} = 0$ a.s., we deduce that $\lim_{n\to\infty}\Delta_n^{(\theta)} = 0$ and $\lim_{n\to\infty}\Delta_n^{(\mu)} = 0$ a.s. Set $T^*$ and $M^*$ such that $\frac{\mathbb{1}_{b=1}}{2\beta_0} < T^* < T$ and $0 < M^* < M$, respectively; we can then deduce from (24) that, for $n$ large enough,

$$
\begin{aligned}
(26) \quad \|\Delta_{n+1}^{(\theta)}\|_T \leq{} & (1 - \beta_n T^*)\|\Delta_n^{(\theta)}\|_T \\
& + \beta_n O[(\beta_n\|L_n^{(\theta)}\| + \beta_n\|R_n^{(\theta)}\| + \|r_n^{(\theta)}\| + \|r_n^{(\mu)}\|)] \\
& + \beta_n C_1(\|L_n^{(\theta)}\|^2 + \|R_n^{(\theta)}\|^2 + \|L_n^{(\mu)}\|^2 + \|R_n^{(\mu)}\|^2 + \|\Delta_n^{(\mu)}\|_M^2)
\end{aligned}
$$

and from (25), that there exists $C_2' > 0$ such that, for $n$ large enough,

$$
\begin{aligned}
(27) \quad \|\Delta_{n+1}^{(\mu)}\|_M \leq{} & (1 - \gamma_n M^*)\|\Delta_n^{(\mu)}\|_M \\
& + \gamma_n[O(\gamma_n\|L_n^{(\mu)}\| + \gamma_n\|R_n^{(\mu)}\| + \|r_n^{(\mu)}\|)] \\
& + \gamma_n M^* C_2'(\|L_n^{(\theta)}\|^2 + \|R_n^{(\theta)}\|^2 + \|L_n^{(\mu)}\|^2 \\
& \quad + \|R_n^{(\mu)}\|^2 + \|\Delta_n^{(\theta)}\|_T).
\end{aligned}
$$

REMARK 5.   Let us note here that classical techniques allow to deduce from (26) that if the sequence

$$(\beta_n\|L_n^{(\theta)}\| + \beta_n\|R_n^{(\theta)}\| + \|r_n^{(\theta)}\| + \|r_n^{(\mu)}\| + \|L_n^{(\theta)}\|^2 + \|R_n^{(\theta)}\|^2$$
$$+ \|L_n^{(\mu)}\|^2 + \|R_n^{(\mu)}\|^2 + \|\Delta_n^{(\mu)}\|_M^2)$$

is bounded above by a suitable sequence $(w'_n)$, then $\|\Delta_n^{(\theta)}\|_T$ can also be bounded above by $(w'_n)$. However, since the first upper bound of $\|\Delta_n^{(\mu)}\|_M$, which will be available in the body of the proof of Lemmas 3 and 4 (see Section 2.3.2) is $\|\Delta_n^{(\mu)}\|_M = O(1)$, inequality (26) leads only to $\|\Delta_n^{(\theta)}\|_T = O(1)$ (which has already been proved). The idea to deduce from (26) a better upper bound for $\|\Delta_n^{(\theta)}\|_T$ is to resort to the averaging effect; for that, we need to substitute $\gamma_n^{-1}[\|\Delta_n^{(\mu)}\|_M - \|\Delta_{n+1}^{(\mu)}\|_M]$ for $\|\Delta_n^{(\mu)}\|_M^2$ in (26) [see (28) and Remark 6 below].

Inequality (27) allows to write

$$\|\Delta_n^{(\mu)}\|_M \leq \frac{1}{\gamma_n M^*}[\|\Delta_n^{(\mu)}\|_M - \|\Delta_{n+1}^{(\mu)}\|_M] + O(\gamma_n\|L_n^{(\mu)}\| + \gamma_n\|R_n^{(\mu)}\| + \|r_n^{(\mu)}\|)$$
$$+ C_2'(\|L_n^{(\theta)}\|^2 + \|R_n^{(\theta)}\|^2 + \|L_n^{(\mu)}\|^2 + \|R_n^{(\mu)}\|^2 + \|\Delta_n^{(\theta)}\|_T).$$

Set $\varepsilon > 0$ such that $\frac{\mathbb{1}_{b=1}}{2\beta_0} < T^* - \varepsilon C_2'$; since $\lim_{n\to\infty}\Delta_n^{(\mu)} = 0$, we deduce from (26) that, for $n$ large enough,

$$\|\Delta_{n+1}^{(\theta)}\|_T \leq (1 - \beta_n T^*)\|\Delta_n^{(\theta)}\|_T$$
$$+ \beta_n[O(\beta_n\|L_n^{(\theta)}\| + \beta_n\|R_n^{(\theta)}\| + \|r_n^{(\theta)}\| + \|r_n^{(\mu)}\|)]$$
$$+ \beta_n C_1(\|L_n^{(\theta)}\|^2 + \|R_n^{(\theta)}\|^2 + \|L_n^{(\mu)}\|^2 + \|R_n^{(\mu)}\|^2) + \beta_n\varepsilon\|\Delta_n^{(\mu)}\|_M$$
$$\leq (1 - \beta_n T^*)\|\Delta_n^{(\theta)}\|_T$$
$$+ \beta_n[O(\beta_n\|L_n^{(\theta)}\| + \beta_n\|R_n^{(\theta)}\| + \|r_n^{(\theta)}\| + \|r_n^{(\mu)}\|)]$$
$$+ \beta_n C_1(\|L_n^{(\theta)}\|^2 + \|R_n^{(\theta)}\|^2 + \|L_n^{(\mu)}\|^2 + \|R_n^{(\mu)}\|^2)$$
$$+ \frac{\beta_n\varepsilon}{\gamma_n M^*}[\|\Delta_n^{(\mu)}\|_M - \|\Delta_{n+1}^{(\mu)}\|_M]$$
$$+ \beta_n[O(\gamma_n\|L_n^{(\mu)}\| + \gamma_n\|R_n^{(\mu)}\| + \|r_n^{(\mu)}\|)]$$
$$+ \beta_n\varepsilon C_2'(\|L_n^{(\theta)}\|^2 + \|R_n^{(\theta)}\|^2 + \|L_n^{(\mu)}\|^2 + \|R_n^{(\mu)}\|^2 + \|\Delta_n^{(\theta)}\|_T).$$

Setting $T^{**}$ such that $\frac{\bar{q}_{b=1}}{2\beta_0} < T^{**} < T^* - \varepsilon C_2'$, we obtain, for $n$ large enough,

$$
\begin{aligned}
(28) \quad \|\Delta_{n+1}^{(\theta)}\|_T \leq &(1 - \beta_n T^{**})\|\Delta_n^{(\theta)}\|_T \\
&+ \beta_n[O(\beta_n\|L_n^{(\theta)}\| + \beta_n\|R_n^{(\theta)}\| + \|r_n^{(\theta)}\| + \|r_n^{(\mu)}\| \\
&\qquad\qquad + \gamma_n\|L_n^{(\mu)}\| + \gamma_n\|R_n^{(\mu)}\|)] \\
&+ \beta_n[O(\|L_n^{(\theta)}\|^2 + \|R_n^{(\theta)}\|^2 + \|L_n^{(\mu)}\|^2 + \|R_n^{(\mu)}\|^2)] \\
&+ \frac{\beta_n\varepsilon}{\gamma_n M^*}[\|\Delta_n^{(\mu)}\|_M - \|\Delta_{n+1}^{(\mu)}\|_M].
\end{aligned}
$$

Classical computations (see Section A.4.2) then allow to deduce from (28) that

$$
(29) \quad \|\Delta_n^{(\theta)}\|_T = O(\beta_n^2\gamma_n^{-2}w_n^2 + \beta_n\gamma_n^{-1}\delta_n^{(\mu)}) + o(\sqrt{\beta_n}).
$$

REMARK 6. Let us point out the averaging effect here again: the term

$$
\gamma_n^{-1}[\|\Delta_n^{(\mu)}\|_M - \|\Delta_{n+1}^{(\mu)}\|_M]
$$

present in (28) leads to the bounding term $\beta_n\gamma_n^{-1}\delta_n^{(\mu)}$ in (29), although the term $\|\Delta_n^{(\mu)}\|_M$ itself is bounded only by $\delta_n^{(\mu)}$.

To conclude the proof of Lemma 6, we substitute the upper bound obtained in (29) for $\|\Delta_n^{(\theta)}\|_T$ in (27) and, via classical computations (see Section A.4.3), establish that

$$
(30) \quad \|\Delta_n^{(\mu)}\|_M = O(\beta_n^2\gamma_n^{-2}w_n^2 + \beta_n\gamma_n^{-1}\delta_n^{(\mu)}) + o(\sqrt{\beta_n}).
$$

Lemma 6 then straightforwardly follows from the equivalence of the finite-dimensional vector norms.

2.3.2. *Body of the proof of Lemmas 3 and 4.* Let $(w_n)$ be a sequence satisfying Conditions (C) and (C′), and such that $\|\mu_n\| = O(w_n)$ a.s. In the proof of Lemma 6, we have seen that $\lim_{n\to\infty}\Delta_n^{(\mu)} = 0$ a.s. We can thus apply Lemma 6 with $(\delta_n^{(\mu)}) \equiv 1$, which ensures that

$$
\|\Delta_n^{(\theta)}\| + \|\Delta_n^{(\mu)}\| = O(\beta_n^2\gamma_n^{-2}w_n^2 + \beta_n\gamma_n^{-1}) + o(\sqrt{\beta_n}) \qquad \text{a.s.}
$$

Now, let $k$ be a positive integer, and assume that

$$
\|\Delta_n^{(\theta)}\| + \|\Delta_n^{(\mu)}\| = O(\beta_n^2\gamma_n^{-2}w_n^2 + [\beta_n\gamma_n^{-1}]^k) + o(\sqrt{\beta_n}) \qquad \text{a.s.}
$$

Since $(w_n)$ satisfies Conditions (C) and (C′), the sequence $(\delta_n^{(\mu)}) \equiv (\beta_n^2 \gamma_n^{-2} w_n^2 + [\beta_n \gamma_n^{-1}]^k + \beta_n^{1/2})$ also satisfies Conditions (C) and (C′); it follows from the application of Lemma 6 that

$$\|\Delta_n^{(\theta)}\| + \|\Delta_n^{(\mu)}\| = O(\beta_n^2 \gamma_n^{-2} w_n^2 + [\beta_n \gamma_n^{-1}]^{k+1}) + o(\sqrt{\beta_n}) \qquad \text{a.s.}$$

We have thus proved by induction that, for all integers $j$,

$$\|\Delta_n^{(\theta)}\| + \|\Delta_n^{(\mu)}\| = O(\beta_n^2 \gamma_n^{-2} w_n^2 + [\beta_n \gamma_n^{-1}]^j) + o(\sqrt{\beta_n}) \qquad \text{a.s.}$$

Since Assumption (A3) ensures the existence of $j_0$ such that $[\beta_n \gamma_n^{-1}]^{j_0} = o(\beta_n^{1/2})$, we have proved that, for any sequence $(w_n)$ satisfying Conditions (C) and (C′) and such that $\|\mu_n\| = O(w_n)$, we have

$$\tag{31} \|\Delta_n^{(\theta)}\| + \|\Delta_n^{(\mu)}\| = O(\beta_n^2 \gamma_n^{-2} w_n^2) + o(\sqrt{\beta_n}) \quad \text{a.s.}$$

Set $k \geq 0$, and assume that

$$\tag{32} \|\mu_n\| = O(\sqrt{\gamma_n \log s_n} + [\beta_n \gamma_n^{-1}]^k) \qquad \text{a.s.}$$

Since the sequence $(\sqrt{\gamma_n \log s_n} + [\beta_n \gamma_n^{-1}]^k)$ satisfies Conditions (C) and (C′), the application of Lemma 2, and of Lemma 5 and (31) with $(w_n) \equiv (\sqrt{\gamma_n \log s_n} + [\beta_n \gamma_n^{-1}]^k)$ ensures that

$$
\begin{aligned}
\|\mu_n\| &= O(\|L_n^{(\mu)}\| + \|R_n^{(\mu)}\| + \|\Delta_n^{(\mu)}\|) \\
&= O(\sqrt{\gamma_n \log s_n} \\
&\qquad + [(\beta_n \gamma_n^{-1})^{k+1} + \beta_n \gamma_n^{-1} \sqrt{\gamma_n \log s_n} + \sqrt{\beta_n \log u_n}] + (\beta_n \gamma_n^{-1})^{2k+2}) \\
&\qquad + o(\sqrt{\beta_n}) \qquad\qquad\qquad \text{a.s.} \\
&= O(\sqrt{\gamma_n \log s_n} + [\beta_n \gamma_n^{-1}]^{k+1}) \qquad \text{a.s. [in view of (A3)].}
\end{aligned}
$$

Now, in view of assumption (A1), we have $\|\mu_n\| = o(1)$ a.s., so that (32) is satisfied for $k = 0$. We have thus proved by induction that (32) holds for all $k \geq 0$. Since (A3) ensures the existence of $k_0$ such that $[\beta_n \gamma_n^{-1}]^{k_0} = o(\sqrt{\gamma_n \log s_n})$, it follows that $\|\mu_n\| = O(\sqrt{\gamma_n \log s_n})$ a.s.

REMARK 7. This latter upper bound of $(\mu_n)$ proves the second assertion of Proposition 1.

To conclude the proof of Lemma 3, we now apply Lemma 5 with $(w_n) \equiv (\sqrt{\gamma_n \log s_n})$:

- For all $s \in ]1/2, \beta_0 \Lambda^{(H)}[$, we have

$$
\begin{aligned}
\|R_n^{(\theta)}\| &= O(\beta_n \gamma_n^{-1/2} \sqrt{\log s_n} + n^{-s}) \qquad \text{a.s.} \\
&= O(n^{-(b-a/2)} \sqrt{\log n} + n^{-s}) \qquad \text{a.s.,}
\end{aligned}
$$

with, in view of (A3), $b - a/2 > b/2$; the first part of Lemma 3 follows.

• We have

$$\|R_n^{(\mu)}\| = O([\beta_n \gamma_n^{-1}]^{1/2}\sqrt{\beta_n \log s_n} + \sqrt{\beta_n \log u_n}) \qquad \text{a.s.},$$

which, in view of (A3), gives the second part of Lemma 3.

To conclude the proof of Lemma 4, we apply (31) with $(w_n) \equiv (\sqrt{\gamma_n \log s_n})$, which gives

$$\|\Delta_n^{(\theta)}\| + \|\Delta_n^{(\mu)}\| = O([\beta_n \gamma_n^{-1}][\beta_n \log s_n]) + o(\sqrt{\beta_n}) \qquad \text{a.s.}$$

In view of (A3), Lemma 4 follows.

## 3. The averaging principle in the context of two-time-scale stochastic approximation algorithms.

3.1. *Asymptotic efficiency of two-time-scale stochastic approximation algorithms.* The averaging principle has been introduced simultaneously by Ruppert [26] and Polyak [24] in the framework of single-time-scale stochastic approximation algorithms, and their pioneer work has been widely discussed and extended in this context (see, among many others, Yin [27], Delyon and Juditsky [6], Polyak and Juditsky [25], Kushner and Yang [16], Dippon and Renz [7, 8], Duflo [9], Kushner and Yin [17] and Pelletier [23]). Let us recall that the foundations of this principle are the following: (i) there exists an algorithm which converges with the optimal rate; however, in general, this "optimal algorithm" cannot be used because it depends on an unknown parameter; (ii) taking a suitable average of a slowly converging algorithm leads to an "averaged algorithm," which has the same asymptotic behavior as the "optimal algorithm."

To introduce the averaging principle in the context of two-time-scale stochastic approximation algorithms, we first need to define the notion of asymptotic efficiency in this framework, that is, to find out what the optimal convergence rate of the two-time-scale algorithms is. For that purpose, we follow the approach employed in the framework of the single-time-scale stochastic approximation algorithms, and consider the class of matricial and two-time-scale algorithms defined as

$$(33) \qquad \theta_{n+1} = \theta_n + \frac{A_\theta}{n} X_{n+1},$$

$$(34) \qquad \mu_{n+1} = \mu_n + \frac{A_\mu}{n^a} Y_{n+1},$$

where $a \in \,]1/2, 1[$, and where $A_\theta$ (resp. $A_\mu$) is a $d \times d$ (resp. $d' \times d'$) nonsingular matrix such that the matrix $A_\theta H + I/2$ (resp. $A_\mu Q_{22}$) is attractive [recall that $H$ and $Q_{22}$ are defined in (A2)]. Following the computations made in the beginning of Section 2.3, and setting $(\beta_n) \equiv (n^{-1})$ and $(\gamma_n) \equiv (n^{-a})$, we rewrite (33)–(34) as

$$(35) \qquad \theta_{n+1} = \theta_n + A_\theta \beta_n (Q_{11}\theta_n + Q_{12}\mu_n + \rho_n^{(\theta)} + V_{n+1}),$$

$$(36) \qquad \mu_{n+1} = \mu_n + A_\mu \gamma_n (Q_{21}\theta_n + Q_{22}\mu_n + \rho_n^{(\mu)} + W_{n+1}).$$

From (36), we get

$$\mu_n = Q_{22}^{-1} A_\mu^{-1} \gamma_n^{-1} (\mu_{n+1} - \mu_n) - Q_{22}^{-1} Q_{21}\theta_n - Q_{22}^{-1}\rho_n^{(\mu)} - Q_{22}^{-1} W_{n+1},$$

which, reintroduced in (35), gives

$$\theta_{n+1} = \theta_n + \beta_n (A_\theta H)\theta_n + \beta_n (A_\theta Q_{12} Q_{22}^{-1} A_\mu^{-1}) \gamma_n^{-1}[\mu_{n+1} - \mu_n]$$
$$+ \beta_n A_\theta (V_{n+1} - Q_{12} Q_{22}^{-1} W_{n+1}) + \beta_n A_\theta (\rho_n^{(\theta)} - Q_{12} Q_{22}^{-1} \rho_n^{(\mu)}).$$

Following the proof of Theorem 1, we obtain

$$\sqrt{n}(\theta_n - \theta) \xrightarrow{\mathcal{D}} \mathcal{N}(0, \Sigma_\theta(A_\theta)),$$

where $\Sigma_\theta(A_\theta)$ is the solution of the Lyapounov equation

$$\left[ A_\theta H + \frac{I}{2} \right] \Sigma_\theta(A_\theta) + \Sigma_\theta(A_\theta) \left[ H^T A_\theta^T + \frac{I}{2} \right] = -A_\theta \Gamma_\theta A_\theta^T$$

[$\Gamma_\theta$ being defined in (7)]. Classical computations (see, e.g., [9], page 166) ensure that the optimal choice of $A_\theta$ in (33) is $A_\theta = -H^{-1}$, which leads to the optimal asymptotic covariance matrix $\Sigma_\theta(A_\theta) = H^{-1}\Gamma_\theta[H^{-1}]^T$, and to the following CLT for $\theta_n$:

$$(37) \qquad \sqrt{n}(\theta_n - \theta^*) \xrightarrow{\mathcal{D}} \mathcal{N}(0, H^{-1}\Gamma_\theta[H^{-1}]^T).$$

Therefore, one of the conditions we shall require to say that a general two-time-scale stochastic approximation algorithm of the type (1)–(2) is asymptotically efficient is that its fastest component $\theta_n$ satisfies the CLT (37).

Now, the idea to find out the optimal weak convergence rate for the slowest component $\mu_n$ of the two-time-scale stochastic approximation algorithm (1)–(2) is the following. First, we invert the roles of $\theta_n$ and $\mu_n$, that is, we give to $\mu_n$ the position of the fastest component, and consider the following alternative algorithm to the algorithm (33)–(34):

$$\theta_{n+1} = \theta_n + \frac{A_\theta}{n^a} X_{n+1},$$

$$\mu_{n+1} = \mu_n + \frac{A_\mu}{n} Y_{n+1},$$

where $a \in ]1/2, 1[$. Then, we apply the results previously obtained for the matricial two-time-scale stochastic approximation algorithm (33)–(34). Set

$$
(38) \qquad G = Q_{22} - Q_{21} Q_{11}^{-1} Q_{12},
$$

$$
(39) \begin{aligned}
\Gamma_\mu &= \lim_{n \to \infty} \mathbb{E}([W_{n+1} - Q_{21} Q_{11}^{-1} V_{n+1}][W_{n+1} - Q_{21} Q_{11}^{-1} V_{n+1}]^T | \mathcal{F}_n) \\
&= \Gamma_{22} + Q_{21} Q_{11}^{-1} \Gamma_{11} [Q_{11}^{-1}]^T Q_{21}^T - \Gamma_{21} [Q_{11}^{-1}]^T Q_{21}^T - Q_{21} Q_{11}^{-1} \Gamma_{12},
\end{aligned}
$$

and assume that the matrices $A_\mu G + I/2$ and $A_\theta Q_{11}$ are attractive. Following the proof of (37), we deduce that the optimal choice of $A_\mu$ is $A_\mu = -G^{-1}$, which leads to the optimal covariance matrix $G^{-1} \Gamma_\mu [G^{-1}]^T$ and to the following CLT for $\mu_n$:

$$
\sqrt{n}(\mu_n - \mu^*) \xrightarrow{\mathcal{D}} \mathcal{N}(0, G^{-1} \Gamma_\mu [G^{-1}]^T).
$$

We can now precisely define the notion of asymptotical efficiency for two-time-scale stochastic approximation algorithms.

DEFINITION 3. Let $(\tilde{\theta}_n, \tilde{\mu}_n)$ be given by a two-time-scale stochastic approximation algorithm used for the search of the common zero $(\theta^*, \mu^*)$ of two functions $f$ and $g$. Assume that $f$ and $g$ satisfy assumption (A2)(i), and that the error-contaminated observations $(X_{n+1})$ and $(Y_{n+1})$ of $f(\tilde{\theta}_n, \tilde{\mu}_n)$ and $g(\tilde{\theta}_n, \tilde{\mu}_n)$ satisfy assumption (A4). We say that the two-time-scale stochastic approximation algorithm which defines $(\tilde{\theta}_n, \tilde{\mu}_n)$ is asympotically efficient if the two following properties hold:

(P1) $\qquad\qquad \sqrt{n}(\tilde{\theta}_n - \theta^*) \xrightarrow{\mathcal{D}} \mathcal{N}(0, H^{-1} \Gamma_\theta [H^{-1}]^T),$

(P2) $\qquad\qquad \sqrt{n}(\tilde{\mu}_n - \mu^*) \xrightarrow{\mathcal{D}} \mathcal{N}(0, G^{-1} \Gamma_\mu [G^{-1}]^T),$

where $H$, $\Gamma_\theta$, $G$ and $\Gamma_\mu$ are defined in (6), (7), (38) and (39) respectively.

Let us note that a sequence $(\tilde{\theta}_n, \tilde{\mu}_n)$ satisfying properties (P1) and (P2) can be obtained, under suitable assumptions, by simultaneously running the two following two-time-scale stochastic approximation algorithms:

$$
\tilde{\theta}_{n+1} = \tilde{\theta}_n - \frac{H^{-1}}{n} X_{n+1}^{(1)},
$$

$$
\mu_{n+1} = \mu_n + \frac{1}{n^a} Y_{n+1}^{(1)}
$$

and

$$
\theta_{n+1} = \theta_n + \frac{1}{n^a} X_{n+1}^{(2)},
$$

$$
\tilde{\mu}_{n+1} = \tilde{\mu}_n - \frac{G^{-1}}{n} Y_{n+1}^{(2)},
$$

where $X_{n+1}^{(1)}$, $Y_{n+1}^{(1)}$, $X_{n+1}^{(2)}$ and $Y_{n+1}^{(2)}$ are error-contaminated observations of $f(\tilde{\theta}_n, \mu_n)$, $g(\tilde{\theta}_n, \mu_n)$, $f(\theta_n, \tilde{\mu}_n)$ and $g(\theta_n, \tilde{\mu}_n)$, respectively. However, this procedure has two main drawbacks. The first one (which is minor) is that it doubles the number of necessary observations. The second one (which is much more important) is that, most of the time, this procedure cannot be used, the matrices $H$ and $G$ being usually unknown.

3.2. *Averaging of two-time-scale stochastic approximation algorithms.* We can now introduce the averaged two-time-scale stochastic approximation algorithm. Applying the averaging principle, we first define the slowly converging two-time-scale algorithm. For that purpose, we let the sequence $(\theta_n, \mu_n)$ be still defined by the recursive equations (1)–(2), but, this time, the step sizes $(\beta_n)$ and $(\gamma_n)$ fulfill the following assumption:

(A'3) $(\beta_n) \equiv (\beta_0 n^{-b})$ and $(\gamma_n) \equiv (\gamma_0 n^{-a})$ with $\beta_0 > 0$, $\gamma_0 > 0$, and $\frac{1}{2} < a < b < 1$.

We then define the averages of $\theta_k$ and $\mu_k$ by setting

$$(40) \qquad \overline{\theta}_n = \frac{1}{n} \sum_{k=1}^{n} \theta_k \quad \text{and} \quad \overline{\mu}_n = \frac{1}{n} \sum_{k=1}^{n} \mu_k.$$

To establish the joint weak convergence rate of $(\overline{\theta}_n)$ and $(\overline{\mu}_n)$, we need to strengthen assumption (A4) into the following condition:

(A'4) Assumption (A4) is fulfilled with $\|r_n^{(\theta)}\| + \|r_n^{(\mu)}\| = o(n^{-1/2})$.

Our main result in this section is the following theorem.

THEOREM 2 [Joint weak convergence rate of $(\overline{\theta}_n)$ and $(\overline{\mu}_n)$]. *Let $(\theta_n, \mu_n)$ be defined by the recursive equations (1)–(2), and $(\overline{\theta}_n, \overline{\mu}_n)$ by (40). Under assumptions* (A1), (A2), (A'3) *and* (A'4), *we have*

$$(41) \qquad \sqrt{n} \begin{pmatrix} \overline{\theta}_n - \theta^* \\ \overline{\mu}_n - \mu^* \end{pmatrix} \xrightarrow{\mathcal{D}} \mathcal{N}(0, DP\Gamma P^T D^T),$$

*where $\Gamma$ is defined in* (A4)(ii), *and where*

$$D = \begin{pmatrix} H^{-1} & 0 \\ 0 & G^{-1} \end{pmatrix}, \qquad P = \begin{pmatrix} I & -Q_{12}Q_{22}^{-1} \\ -Q_{21}Q_{11}^{-1} & I \end{pmatrix}.$$

*In particular, the averaged two-time-scale stochastic approximation algorithm $(\overline{\theta}_n, \overline{\mu}_n)$ is asymptotically efficient.*

3.3. *Proof of Theorem* 2. Let us first note that the CLT (41) implies, in particular, that

$$\sqrt{n}(\overline{\theta}_n - \theta^*) \xrightarrow{\mathcal{D}} \mathcal{N}(0, H^{-1}\Gamma_\theta[H^{-1}]^T),$$

$$\sqrt{n}(\overline{\mu}_n - \mu^*) \xrightarrow{\mathcal{D}} \mathcal{N}(0, G^{-1}\Gamma_\mu[G^{-1}]^T),$$

which proves the asymptotic efficiency of the averaged algorithm $(\overline{\theta}_n, \overline{\mu}_n)$. We now prove (41).

We assume again, without loss of generality, that $\theta^* = 0$ and $\mu^* = 0$. In the beginning of Section 2.3 we have seen that [see (13)]:

$$\theta_{n+1} = \theta_n + \beta_n H \theta_n + \beta_n Q_{12} Q_{22}^{-1} \gamma_n^{-1} [\mu_{n+1} - \mu_n] + \beta_n (V_{n+1} - Q_{12} Q_{22}^{-1} W_{n+1})$$
$$+ \beta_n ([\rho_n^{(\theta)} + r_n^{(\theta)}] - Q_{12} Q_{22}^{-1} [\rho_n^{(\mu)} + r_n^{(\mu)}]).$$

We can thus write

$$\theta_n = -H^{-1}(V_{n+1} - Q_{12} Q_{22}^{-1} W_{n+1}) + H^{-1}\beta_n^{-1}[\theta_{n+1} - \theta_n]$$
$$- H^{-1} Q_{12} Q_{22}^{-1} \gamma_n^{-1}[\mu_{n+1} - \mu_n]$$
$$- H^{-1}([\rho_n^{(\theta)} + r_n^{(\theta)}] - Q_{12} Q_{22}^{-1}[\rho_n^{(\mu)} + r_n^{(\mu)}]),$$

so that

$$\overline{\theta}_n = H^{-1}\left( -\frac{1}{n}\sum_{k=1}^{n}(V_{k+1} - Q_{12} Q_{22}^{-1} W_{k+1}) + \mathcal{R}_n^{(1)} - \mathcal{R}_n^{(2)} - \mathcal{R}_n^{(3)} - \mathcal{R}_n^{(4)} \right),$$

with

$$\mathcal{R}_n^{(1)} = \frac{1}{n}\sum_{k=1}^{n} \beta_k^{-1}[\theta_{k+1} - \theta_k],$$

$$\mathcal{R}_n^{(2)} = \frac{1}{n}\sum_{k=1}^{n} Q_{12} Q_{22}^{-1} \gamma_k^{-1}[\mu_{k+1} - \mu_k],$$

$$\mathcal{R}_n^{(3)} = \frac{1}{n}\sum_{k=1}^{n}[\rho_k^{(\theta)} - Q_{12} Q_{22}^{-1} \rho_k^{(\mu)}],$$

$$\mathcal{R}_n^{(4)} = \frac{1}{n}\sum_{k=1}^{n}[r_k^{(\theta)} - Q_{12} Q_{22}^{-1} r_k^{(\mu)}].$$

Similarly, we have

$$\overline{\mu}_n = G^{-1}\left( -\frac{1}{n}\sum_{k=1}^{n}(W_{k+1} - Q_{21} Q_{11}^{-1} V_{k+1}) + \mathcal{R}_n^{(5)} - \mathcal{R}_n^{(6)} - \mathcal{R}_n^{(7)} - \mathcal{R}_n^{(8)} \right),$$

with

$$\mathcal{R}_n^{(5)} = \frac{1}{n} \sum_{k=1}^{n} \gamma_k^{-1} [\mu_{k+1} - \mu_k],$$

$$\mathcal{R}_n^{(6)} = \frac{1}{n} \sum_{k=1}^{n} Q_{21} Q_{11}^{-1} \beta_k^{-1} [\theta_{k+1} - \theta_k],$$

$$\mathcal{R}_n^{(7)} = \frac{1}{n} \sum_{k=1}^{n} [\rho_k^{(\mu)} - Q_{21} Q_{11}^{-1} \rho_k^{(\theta)}],$$

$$\mathcal{R}_n^{(8)} = \frac{1}{n} \sum_{k=1}^{n} [r_k^{(\mu)} - Q_{21} Q_{11}^{-1} r_k^{(\theta)}].$$

A straightforward application of Lyapounov's theorem gives the following lemma:

LEMMA 7.

$$\frac{1}{\sqrt{n}} \sum_{k=1}^{n} \begin{pmatrix} V_{k+1} - Q_{12} Q_{22}^{-1} W_{k+1} \\ W_{k+1} - Q_{21} Q_{11}^{-1} V_{k+1} \end{pmatrix} \xrightarrow{\mathcal{D}} \mathcal{N}(0, P\Gamma P^T).$$

The CLT (41) follows thus from the combination of Lemma 7 and of the following lemma.

LEMMA 8. *For $i \in \{1, \ldots, 8\}$, we have*

$$\lim_{n \to \infty} \sqrt{n} \mathcal{R}_n^{(i)} = 0 \qquad a.s.$$

PROOF. The application of Proposition 1 gives

$$\frac{1}{\sqrt{n}} \sum_{k=1}^{n} \beta_k^{-1} [\theta_{k+1} - \theta_k]$$

$$= \frac{1}{\sqrt{n}} \left( \frac{\theta_{n+1}}{\beta_n} - \frac{\theta_1}{\beta_1} + \sum_{k=1}^{n} \left[ \frac{1}{\beta_{k-1}} - \frac{1}{\beta_k} \right] \theta_k \right)$$

$$= O\left( \frac{\|\theta_{n+1}\|}{\sqrt{n}\beta_n} + \frac{1}{\sqrt{n}} + \frac{1}{\sqrt{n}} \sum_{k=1}^{n} k^{b-1} \|\theta_k\| \right)$$

$$= O\left( \frac{\sqrt{\beta_n \log u_n}}{\sqrt{n}\beta_n} + \frac{1}{\sqrt{n}} + \frac{1}{\sqrt{n}} \sum_{k=1}^{n} k^{b-1} \sqrt{\beta_k \log u_k} \right) \qquad \text{a.s.}$$

$$= O(n^{b/2 - 1/2} \log n) \qquad \text{a.s.}$$

Since $b < 1$, it follows that $\lim_{n\to\infty} \sqrt{n}\mathcal{R}_n^{(1)} = 0$ and $\lim_{n\to\infty} \sqrt{n}\mathcal{R}_n^{(6)} = 0$ a.s. In the same way, we have

$$\frac{1}{\sqrt{n}}\sum_{k=1}^{n}\gamma_k^{-1}[\mu_{k+1} - \mu_k]$$

$$= \frac{1}{\sqrt{n}}\left(\frac{\mu_{n+1}}{\gamma_n} - \frac{\mu_1}{\gamma_1} + \sum_{k=1}^{n}\left[\frac{1}{\gamma_{k-1}} - \frac{1}{\gamma_k}\right]\mu_k\right)$$

$$= O\left(\frac{\sqrt{\log s_n}}{\sqrt{n\gamma_n}} + \frac{1}{\sqrt{n}} + \frac{1}{\sqrt{n}}\sum_{k=1}^{n}k^{a-1}\sqrt{\gamma_k\log s_k}\right) \qquad \text{a.s.}$$

$$= O(n^{-(1-a)/2}\log n) \qquad\qquad\qquad \text{a.s.}$$

Since $a < 1$, it follows that $\lim_{n\to\infty} \sqrt{n}\mathcal{R}_n^{(2)} = 0$ and $\lim_{n\to\infty} \sqrt{n}\mathcal{R}_n^{(5)} = 0$ a.s. Now, we note that

$$\frac{1}{\sqrt{n}}\sum_{k=1}^{n}(\|\theta_k\|^2 + \|\mu_k\|^2) = O\left(\frac{1}{\sqrt{n}}\sum_{k=1}^{n}\gamma_k\log s_k\right) \qquad \text{a.s.}$$

$$= O(n^{1/2-a}\log n) \qquad \text{a.s.}$$

Since $a > 1/2$ and in view of (12), it follows that $\lim_{n\to\infty}\sqrt{n}\mathcal{R}_n^{(3)} = 0$ and $\lim_{n\to\infty}\sqrt{n}\mathcal{R}_n^{(7)} = 0$ a.s. Finally, assumption (A'4) ensures that $\lim_{n\to\infty}\sqrt{n}\times\mathcal{R}_n^{(4)} = 0$ and $\lim_{n\to\infty}\sqrt{n}\mathcal{R}_n^{(8)} = 0$ a.s. $\quad\square$

# APPENDIX

## A.1. Two technical lemmas.

LEMMA 9. *Let* $(x_n)$ *be a sequence of positive real numbers, let* $(u_n)$ *be an* $\mathbb{R}^d$-*valued random sequence such that* $\|u_n\| = O(x_n)$ *a.s., set* $T > 0$,

$$Z_n^{(1)} = e^{-u_n T}\sum_{k=1}^{n}e^{u_k T}\beta_k x_k \quad \text{and} \quad Z_n^{(2)} = e^{u_n H}\sum_{k=1}^{n}e^{-u_k H}\beta_k u_k.$$

*Let* $(w_n)$ *be a nonrandom sequence satisfying Condition* (C).

1. *For all* $T' \in \,]0, T[$, *we have*

$$|Z_n^{(1)}| = \begin{cases} O(e^{-u_n T'}\mathbb{1}_{b=1} + w_n), & \text{if } x_n = O(w_n), \\ o(e^{-u_n T'}\mathbb{1}_{b=1} + w_n), & \text{if } x_n = o(w_n). \end{cases}$$

2. *For all* $T' \in \,]0, \Lambda^{(H)}[$, *we have*

$$\|Z_n^{(2)}\| = \begin{cases} O(e^{-u_n T'}\mathbb{1}_{b=1} + w_n), & \text{if } x_n = O(w_n), \\ o(e^{-u_n T'}\mathbb{1}_{b=1} + w_n), & \text{if } x_n = o(w_n). \end{cases}$$

LEMMA 10. *Let $(x_n)$ be a sequence of positive real numbers, let $(u_n)$ be an $\mathbb{R}^d$-valued random sequence such that $\|u_n\| = O(x_n)$ a.s., set $T > 0$,*

$$Z_n^{(1)} = e^{-s_n T} \sum_{k=1}^n e^{s_k T} \gamma_k x_k \quad and \quad Z_n^{(2)} = e^{s_n Q_{22}} \sum_{k=1}^n e^{-s_k Q_{22}} \gamma_k u_k.$$

*Let $(w_n)$ be a nonrandom sequence satisfying Condition (C′). We have*

$$|Z_n^{(1)}| + \|Z_n^{(2)}\| = \begin{cases} O(w_n), & \text{if } x_n = O(w_n), \\ o(w_n), & \text{if } x_n = o(w_n). \end{cases}$$

PROOF OF LEMMA 9.  We first establish the upper bound of $(Z_n^{(1)})$.

• Consider the case $b = 1$, that is, $(\beta_n) \equiv (\beta_0 n^{-1})$. In the case $x_n = o(w_n)$, we have

$$\begin{aligned}
|Z_n^{(1)}| &= O\left(n^{-\beta_0 T} \sum_{k=1}^n k^{\beta_0 T - 1} x_k\right) \\
&= o\left(n^{-\beta_0 T} \sum_{k=1}^n k^{\beta_0 T - 1 - \omega} \mathcal{L}(k)\right) \\
&= o(n^{-\beta_0 T}[\log n + n^{\beta_0 T - \omega}] \mathcal{L}(n)) \\
&= o(n^{-\beta_0 T} \mathcal{L}(n) \log n + w_n).
\end{aligned}$$

Since $\mathcal{L}$ is a slowly varying function, it follows that, for all $T' \in ]0, T[$,

$$\begin{aligned}
|Z_n^{(1)}| &= o(n^{-\beta_0 T'} + w_n) \\
&= o(e^{-u_n T'} + w_n).
\end{aligned}$$

In the case $x_n = O(w_n)$, the upper bound of $(Z_n^{(1)})$ is obtained by replacing $o(\cdot)$ by $O(\cdot)$ in the previous equations.

• Consider the case $b < 1$. We note that the sequence $(Z_n^{(1)})$ satisfies the recursive equation

$$Z_n^{(1)} = e^{\beta_n T} Z_{n-1}^{(1)} + \beta_n x_n,$$

so that we can write

$$\begin{aligned}
w_n^{-1} Z_n^{(1)} &= e^{-\beta_n T} \left[\frac{w_{n-1}}{w_n}\right] (w_{n-1}^{-1} Z_{n-1}^{(1)}) + \beta_n w_n^{-1} x_n \\
&= [1 - \beta_n T + O(\beta_n^2)][1 + o(\beta_n)](w_{n-1}^{-1} Z_{n-1}^{(1)}) + \beta_n w_n^{-1} x_n \\
&= [1 - \beta_n T + o(\beta_n)](w_{n-1}^{-1} Z_{n-1}^{(1)}) + \beta_n w_n^{-1} x_n.
\end{aligned}$$

Now, set $T' \in ]0, T[$; for $n$ large enough, we get

$$|w_n^{-1} Z_n^{(1)}| \leq (1 - \beta_n T') |w_{n-1}^{-1} Z_{n-1}^{(1)}| + \beta_n w_n^{-1} x_n,$$

and the application of Lemma 4.I.2 of [9] ensures that if $x_n = O(w_n)$, then the sequence $(w_{n-1}^{-1} Z_{n-1}^{(1)})$ is bounded, that is, $|Z_n^{(1)}| = O(w_n)$; if $x_n = o(w_n)$, then the sequence $(w_{n-1}^{-1} Z_{n-1}^{(1)})$ goes to zero, that is, $|Z_n^{(1)}| = o(w_n)$.

We now establish the upper bound of $(Z_n^{(2)})$. Let $\|\| \cdot \|\|$ denote the matrix norm associated with the Euclidean vector norm. We have

$$\|Z_n^{(2)}\| \leq \sum_{k=1}^{n} \|\|e^{(u_n - u_k)H}\|\| \beta_k \|u_k\|,$$

and the application of Proposition 3.I.2 of [9] ensures that, for all $T \in ]0, \Lambda^{(H)}[$,

$$\|Z_n^{(2)}\| \leq \sum_{k=1}^{n} e^{-(u_n - u_k)T} \beta_k w_k \qquad \text{a.s.}$$

The upper bound of $(Z_n^{(2)})$ then follows straightforwardly from the one obtained for $(Z_n^{(1)})$. $\quad\square$

PROOF OF LEMMA 10. The proof is straightforward by following the proof of Lemma 9 in the case $b < 1$. $\quad\square$

**A.2. Proof of Lemma 1.** Set

$$M_j^{(n)} = \begin{pmatrix} \sqrt{\beta_n^{-1}} e^{u_n H} & 0 \\ 0 & \sqrt{\gamma_n^{-1}} e^{s_n Q_{22}} \end{pmatrix}$$
$$\times \sum_{k=1}^{j} \begin{pmatrix} e^{-u_k H} \beta_k (V_{k+1} - Q_{12} Q_{22}^{-1} W_{k+1}) \\ e^{-s_k Q_{22}} \gamma_k W_{k+1} \end{pmatrix}.$$

For each $n$, $M^{(n)} = (M_j^{(n)})_{j \geq 1}$ is a martingale whose increasing process satisfies

$$\langle M \rangle_n^{(n)} = \begin{pmatrix} A_{1,n} & A_{2,n} \\ A_{2,n}^T & A_{4,n} \end{pmatrix},$$

with, in view of assumption (A4),

$$A_{1,n} = \beta_n^{-1} e^{u_n H} \left\{ \sum_{k=1}^{n} \beta_k^2 e^{-u_k H} \Gamma_\theta e^{-u_k H^T} \right\} e^{u_n H^T},$$

$$A_{2,n} = \sqrt{\beta_n^{-1} \gamma_n^{-1}} e^{u_n H} e^{-s_n} \left\{ \sum_{k=1}^{n} \beta_k \gamma_k e^{-u_k H} \Gamma_{1,2} e^{-s_k Q_{22}^T} \right\} e^{s_n Q_{22}^T},$$

$$A_{4,n} = \gamma_n^{-1} e^{s_n Q_{22}} \left\{ \sum_{k=1}^{n} \gamma_k^2 e^{-s_k Q_{22}} \Gamma_\mu e^{-s_k Q_{22}^T} \right\} e^{s_n Q_{22}^T}.$$

The application of Lemma 4 in [19] ensures that

$$\lim_{n\to\infty} A_{1,n} = \Sigma_\theta \quad \text{and} \quad \lim_{n\to\infty} A_{4,n} = \Sigma_\mu.$$

Moreover, we note that

$$\|A_{2,n}\| = O\left(\sqrt{\beta_n^{-1}\gamma_n^{-1}} \sum_{k=1}^{n} \beta_k \gamma_k \|e^{(u_n-u_k)H}\| \|\| e^{(s_n-s_k)Q_{22}^T}\|\right).$$

Set $T \in ]0, \Lambda^{(H)}[$ and $T' \in ]0, \Lambda^{(Q_{22})}[$; the application of Proposition 3.I.2 in [9] ensures that

$$\|A_{2,n}\| = O\left(\sqrt{\beta_n^{-1}\gamma_n^{-1}} \sum_{k=1}^{n} \beta_k \gamma_k e^{-T(u_n-u_k)} e^{-T'(s_n-s_k)}\right)$$

$$= O\left(\sqrt{\beta_n^{-1}\gamma_n^{-1}} \sum_{k=1}^{n} \beta_k \gamma_k e^{-T'(s_n-s_k)}\right),$$

and the application of Lemma 10 gives

$$\|A_{2,n}\| = O(\sqrt{\beta_n^{-1}\gamma_n^{-1}} \beta_n) = O(\sqrt{\beta_n \gamma_n^{-1}}).$$

In view of (A3), it follows that $\lim_{n\to\infty} A_{2,n} = 0$, and we thus obtain

$$\lim_{n\to\infty} \langle M \rangle_n^{(n)} = \begin{pmatrix} \Sigma_\theta & 0 \\ 0 & \Sigma_\mu \end{pmatrix} \quad \text{a.s.}$$

Now, set $T \in ]\mathbb{1}_{b=1}/(2\beta_0), \Lambda^{(H)}[$ and $T' \in ]0, \Lambda^{(Q_{22})}[$; in view of assumption (A4), we have

$$\sum_{k=1}^{n} \mathbb{E}[\|M_k^{(n)} - M_{k-1}^{(n)}\|^m | \mathcal{F}_{k-1}]$$

$$= O\left(\sum_{k=1}^{n} (\beta_n^{-m/2} \|\beta_k e^{(u_n-u_k)H}\|^m) + \sum_{k=1}^{n} \gamma_n^{-m/2} \|\gamma_k e^{(s_n-s_k)Q_{22}}\|^m\right) \quad \text{a.s.}$$

$$= O\left(\sum_{k=1}^{n} (\beta_n^{-m/2} \beta_k^m e^{-mT(u_n-u_k)}) + \sum_{k=1}^{n} \gamma_n^{-m/2} \gamma_k^m e^{-mT'(s_n-s_k)}\right) \quad \text{a.s.},$$

where the latter upper bound follows from the application of Proposition 3.I.2 in [9]. The application of Lemmas 9 and 10 then ensures that, for all $T^* \in ]\mathbb{1}_{b=1}/(2\beta_0), T[$,

$$\sum_{k=1}^{n} \mathbb{E}[\|M_k^{(n)} - M_{k-1}^{(n)}\|^m | \mathcal{F}_{k-1}]$$

$$= O(\beta_n^{-m/2}[e^{-mT^* u_n} \mathbb{1}_{b=1} + \beta_n^{m-1}] + \gamma_n^{-m/2+(m-1)}) \quad \text{a.s.}$$

$$= O(n^{m/2-\beta_0 mT^*} \mathbb{1}_{b=1} + \beta_n^{m/2-1} + \gamma_n^{m/2-1}) \quad \text{a.s.},$$

so that it be comes

$$\lim_{n\to\infty}\sum_{k=1}^{n}\mathbb{E}[\|M_k^{(n)}-M_{k-1}^{(n)}\|^m|\mathcal{F}_{k-1}]=0 \qquad \text{a.s.}$$

The application of Lyapounov's theorem then gives

$$M_n^{(n)}=\begin{pmatrix}\sqrt{\beta_n^{-1}}L_{n+1}^{(\theta)}\\\sqrt{\gamma_n^{-1}}L_{n+1}^{(\mu)}\end{pmatrix}\xrightarrow{\mathcal{D}}\mathcal{N}\left(0,\begin{pmatrix}\Sigma_\theta & 0\\0 & \Sigma_\mu\end{pmatrix}\right),$$

which concludes the proof of Lemma 1.

**A.3. Proof of Lemma 5.** We first note that, in view of (15), we have

$$R_{n+1}^{(\theta)}=\beta_n\gamma_n^{-1}Q_{12}Q_{22}^{-1}\mu_{n+1}$$

$$+e^{u_nH}\sum_{k=2}^{n}[e^{-u_{k-1}H}\beta_{k-1}\gamma_{k-1}^{-1}-e^{-u_kH}\beta_k\gamma_k^{-1}]Q_{12}Q_{22}^{-1}\mu_k$$

$$-e^{u_nH}e^{-u_1H}\beta_1\gamma_1^{-1}Q_{12}Q_{22}^{-1}\mu_1$$

$$=\beta_n\gamma_n^{-1}Q_{12}Q_{22}^{-1}\mu_{n+1}$$

$$+e^{u_nH}\sum_{k=2}^{n}e^{-u_kH}\beta_kU_k-e^{u_nH}e^{-u_1H}\beta_1\gamma_1^{-1}Q_{12}Q_{22}^{-1}\mu_1,$$

where

$$U_k=\gamma_k^{-1}[e^{-\beta_kH}\beta_{k-1}\beta_k^{-1}\gamma_{k-1}^{-1}\gamma_k-I]Q_{12}Q_{22}^{-1}\mu_k.$$

It follows that

$$\|R_{n+1}^{(\theta)}\|=O\left(\beta_n\gamma_n^{-1}\|\mu_{n+1}\|+\left\|e^{u_nH}\sum_{k=2}^{n}e^{-u_kH}\beta_kU_k\right\|+\|\!|e^{u_nH}|\!\|\right).$$

Note that

$$\|U_n\|=O(\gamma_n^{-1}\beta_n\|\mu_n\|)$$

$$=O(\gamma_n^{-1}\beta_nw_n) \qquad \text{a.s.}$$

Since the sequence $(w_n)$ satisfies Condition (C), the sequence $(\gamma_n^{-1}\beta_nw_n)$ satisfies Condition (C); it follows from the application of Lemma 9 that, for all $t\in]0,\Lambda^{(H)}[$, we have

$$\|R_{n+1}^{(\theta)}\|=O(\beta_n\gamma_n^{-1}w_n+e^{-u_nt}\mathbb{1}_{b=1})+O(\|\!|e^{u_nH}|\!\|) \qquad \text{a.s.}$$

Now, the application of Proposition 3.I.2 in [9] ensures that, for all $t\in]0,\Lambda^{(H)}[$, $\|\!|e^{u_nH}|\!\|=O(e^{-u_nt})$; it follows that, for all $t\in]0,\Lambda^{(H)}[$,

$$\|R_{n+1}^{(\theta)}\|=O(\beta_n\gamma_n^{-1}w_n+e^{-u_nt}) \qquad \text{a.s.}$$

and thus, for all $s \in ]\frac{1}{2}, \beta_0 \Lambda^{(H)}[$, we obtain

$$\|R_{n+1}^{(\theta)}\| = O(\beta_n \gamma_n^{-1} w_n + n^{-s}) \qquad \text{a.s.,}$$

which proves the first part of Lemma 5.

In view of (18), we note that

$$R_{n+1}^{(\mu)} = e^{s_n Q_{22}} \sum_{k=1}^{n} e^{-s_k Q_{22}} \gamma_k \tilde{U}_k,$$

with, by application of Lemma 2 and of the first part of Lemma 5,

$$\begin{aligned}
\|\tilde{U}_n\| &= O(\|L_n^{(\theta)}\| + \|R_n^{(\theta)}\|) \\
&= O(\sqrt{\beta_n \log u_n} + \beta_n \gamma_n^{-1} w_n + n^{-s}) \qquad \text{a.s.} \\
&= O(\sqrt{\beta_n \log u_n} + \beta_n \gamma_n^{-1} w_n) \qquad\qquad \text{a.s.}
\end{aligned}$$

Since the sequence $(w_n)$ satisfies Condition (C'), the sequences $(\sqrt{\beta_n \log u_n})$ and $(\beta_n \gamma_n^{-1} w_n)$ satisfy Condition (C'); the application of Lemma 10 gives

$$\|R_{n+1}^{(\mu)}\| = O(\sqrt{\beta_n \log u_n} + \beta_n \gamma_n^{-1} w_n) \qquad \text{a.s.,}$$

which concludes the proof of Lemma 5.

### A.4. Technical details for the proof of Lemma 6.

A.4.1. *Proof of* (20) *and* (21). Noting that, in view of (14) and (15), we have

$$\begin{aligned}
L_{n+1}^{(\theta)} &= \beta_n(V_{n+1} - Q_{12}Q_{22}^{-1}W_{n+1}) + e^{\beta_n H}L_n^{(\theta)}, \\
R_{n+1}^{(\theta)} &= \beta_n Q_{12}Q_{22}^{-1}\gamma_n^{-1}[\mu_{n+1} - \mu_n] + e^{\beta_n H}R_n^{(\theta)},
\end{aligned}$$

and using (16) and then (13), we write

$$\begin{aligned}
\Delta_{n+1}^{(\theta)} &= \theta_{n+1} - L_{n+1}^{(\theta)} - R_{n+1}^{(\theta)} \\
&= \theta_n + \beta_n H \theta_n + \beta_n([\rho_n^{(\theta)} + r_n^{(\theta)}] - Q_{12}Q_{22}^{-1}[\rho_n^{(\mu)} + r_n^{(\mu)}]) \\
&\quad - e^{\beta_n H}L_n^{(\theta)} - e^{\beta_n H}R_n^{(\theta)} \\
&= (I + \beta_n H)\theta_n - [I + \beta_n H + O(\beta_n^2)][L_n^{(\theta)} + R_n^{(\theta)}] \\
&\quad + \beta_n([\rho_n^{(\theta)} + r_n^{(\theta)}] - Q_{12}Q_{22}^{-1}[\rho_n^{(\mu)} + r_n^{(\mu)}]) \\
&= (I + \beta_n H)\Delta_n^{(\theta)} + O(\beta_n^2)[L_n^{(\theta)} + R_n^{(\theta)}] \\
&\quad + \beta_n([\rho_n^{(\theta)} + r_n^{(\theta)}] - Q_{12}Q_{22}^{-1}[\rho_n^{(\mu)} + r_n^{(\mu)}]),
\end{aligned}$$

which proves (20). Similarly, we note that, in view of (11), (17), (18) and (19), we have

$$
\begin{aligned}
\Delta_{n+1}^{(\mu)} &= \mu_{n+1} - L_{n+1}^{(\mu)} - R_{n+1}^{(\mu)} \\
&= \mu_n + \gamma_n (Q_{21}\theta_n + Q_{22}\mu_n + \rho_n^{(\mu)} + r_n^{(\mu)} + W_{n+1}) - \gamma_n W_{n+1} \\
&\quad - e^{\gamma_n Q_{22}} L_n^{(\mu)} - \gamma_n Q_{21}[L_n^{(\theta)} + R_n^{(\theta)}] - e^{\gamma_n Q_{22}} R_n^{(\mu)}.
\end{aligned}
$$

Using (16), it follows that

$$
\begin{aligned}
\Delta_{n+1}^{(\mu)} &= \mu_n + \gamma_n (Q_{21}\theta_n + Q_{22}\mu_n + \rho_n^{(\mu)} + r_n^{(\mu)}) - e^{\gamma_n Q_{22}} L_n^{(\mu)} \\
&\quad - \gamma_n Q_{21}[\theta_n - \Delta_n^{(\theta)}] - e^{\gamma_n Q_{22}} R_n^{(\mu)} \\
&= (I + \gamma_n Q_{22})\mu_n - [I + \gamma_n Q_{22} + O(\gamma_n^2)][L_n^{(\mu)} + R_n^{(\mu)}] \\
&\quad + \gamma_n[\rho_n^{(\mu)} + r_n^{(\mu)} + Q_{21}\Delta_n^{(\theta)}] \\
&= (I + \gamma_n Q_{22})\Delta_n^{(\mu)} + O(\gamma_n^2)[L_n^{(\mu)} + R_n^{(\mu)}] \\
&\quad + \gamma_n[\rho_n^{(\mu)} + r_n^{(\mu)} + Q_{21}\Delta_n^{(\theta)}],
\end{aligned}
$$

which proves (21).

A.4.2. *Proof of* (29). In view of (28), we have

$$
\|\Delta_{n+1}^{(\theta)}\|_T \le (1 - \beta_n T^{**})\|\Delta_n^{(\theta)}\|_T + \beta_n z_n + \frac{\beta_n \varepsilon}{\gamma_n M^*}[\|\Delta_n^{(\mu)}\|_M - \|\Delta_{n+1}^{(\mu)}\|_M],
$$

where $(z_n)$ is a nonnegative sequence such that

$$
\begin{aligned}
z_n = O(\beta_n \|L_n^{(\theta)}\| &+ \beta_n \|R_n^{(\theta)}\| + \|r_n^{(\theta)}\| + \|r_n^{(\mu)}\| + \gamma_n \|L_n^{(\mu)}\| + \gamma_n \|R_n^{(\mu)}\| \\
&+ \|L_n^{(\theta)}\|^2 + \|R_n^{(\theta)}\|^2 + \|L_n^{(\mu)}\|^2 + \|R_n^{(\mu)}\|^2).
\end{aligned}
$$

For $n \ge 1$, set $\pi_n = \prod_{k=1}^n (1 - \beta_k T^{**})$. We note that

$$
\begin{aligned}
\|\Delta_{n+1}^{(\theta)}\|_T &\le \pi_n \|\Delta_1^{(\theta)}\|_T + \sum_{k=1}^n \frac{\pi_n}{\pi_k} \beta_k \left[ z_k + \frac{\varepsilon}{\gamma_k M^*}(\|\Delta_k^{(\mu)}\|_M - \|\Delta_{k+1}^{(\mu)}\|_M) \right] \\
&\le \pi_n \|\Delta_1^{(\theta)}\|_T + \sum_{k=1}^n \frac{\pi_n}{\pi_k} \beta_k z_k \\
&\quad + \frac{\pi_n \varepsilon}{M^*} \left[ \sum_{k=2}^n \left( \frac{\beta_k}{\pi_k \gamma_k} - \frac{\beta_{k-1}}{\pi_{k-1}\gamma_{k-1}} \right) \|\Delta_k^{(\mu)}\|_M \right. \\
&\quad \left. + \frac{\beta_1}{\pi_1 \gamma_1} \|\Delta_1^{(\mu)}\|_M - \frac{\beta_n}{\pi_n \gamma_n} \|\Delta_{n+1}^{(\mu)}\|_M \right]
\end{aligned}
$$

$$\leq \pi_n \left[ \|\Delta_1^{(\theta)}\|_T + \frac{\beta_1 \varepsilon}{\pi_1 \gamma_1 M^*} \|\Delta_1^{(\mu)}\|_M \right] + \sum_{k=1}^{n} \frac{\pi_n}{\pi_k} \beta_k z_k$$

$$+ \frac{\pi_n \varepsilon}{M^*} \sum_{k=2}^{n} \frac{\beta_k}{\pi_k \gamma_k} \left[ 1 - \frac{\beta_{k-1} \pi_k \gamma_k}{\beta_k \pi_{k-1} \gamma_{k-1}} \right] \|\Delta_k^{(\mu)}\|_M.$$

Since $\pi_k/\pi_{k-1} = 1 - \beta_k T^{**}$ and since, in view of (A3), $\frac{\beta_{k-1} \gamma_k}{\beta_k \gamma_{k-1}} = 1 + O(\beta_k)$, there exists $c > 0$ such that

$$\|\Delta_{n+1}^{(\theta)}\|_T \leq \pi_n \left[ \|\Delta_1^{(\theta)}\|_T + \frac{\beta_1 \varepsilon}{\pi_1 \gamma_1 M^*} \|\Delta_1^{(\mu)}\|_M \right] + \sum_{k=1}^{n} \frac{\pi_n}{\pi_k} \beta_k \left[ z_k + c \frac{\beta_k}{\gamma_k} \|\Delta_k^{(\mu)}\|_M \right].$$

Noting that $\pi_n/\pi_k \leq e^{-T^{**}(u_n - u_k)}$, it follows that

$$\|\Delta_{n+1}^{(\theta)}\|_T = O\left( e^{-T^{**}u_n} + e^{-T^{**}u_n} \sum_{k=1}^{n} e^{T^{**}u_k} \beta_k \left[ z_k + \frac{\beta_k}{\gamma_k} \|\Delta_k^{(\mu)}\|_M \right] \right),$$

which can be rewritten as

$$\|\Delta_{n+1}^{(\theta)}\|_T = O\Bigg( e^{-T^{**}u_n}$$

$$+ e^{-T^{**}u_n} \sum_{k=1}^{n} e^{T^{**}u_k} \beta_k [\beta_k \|L_k^{(\theta)}\| + \beta_k \|R_k^{(\theta)}\| + \|r_k^{(\theta)}\| + \|r_k^{(\mu)}\|]$$

$$+ e^{-T^{**}u_n} \sum_{k=1}^{n} e^{T^{**}u_k} \beta_k [\gamma_k \|L_k^{(\mu)}\| + \gamma_k \|R_k^{(\mu)}\| + \|L_k^{(\theta)}\|^2$$

$$+ \|R_k^{(\theta)}\|^2 + \|L_k^{(\mu)}\|^2 + \|R_k^{(\mu)}\|^2]$$

$$+ e^{-T^{**}u_n} \sum_{k=1}^{n} e^{T^{**}u_k} \frac{\beta_k^2}{\gamma_k} \|\Delta_k^{(\mu)}\|_M \Bigg).$$

Replacing $\|r_k^{(\theta)}\|$ and $\|r_k^{(\mu)}\|$ by their upper bounds given in assumption (A4)(iv), $\|L_k^{(\theta)}\|$ and $\|L_k^{(\mu)}\|$ by their upper bounds given by Lemma 2, $\|R_k^{(\theta)}\|$ and $\|R_k^{(\mu)}\|$ by their upper bounds given by Lemma 5, $\|\Delta_k^{(\mu)}\|$ by $\delta_k^{(\mu)}$, and doing some straightforward simplifications, we obtain

$$\|\Delta_{n+1}^{(\theta)}\|_T = O\Bigg( e^{-T^{**}u_n}$$

$$+ e^{-T^{**}u_n} \sum_{k=1}^{n} e^{T^{**}u_k} \beta_k [\beta_k^2 \gamma_k^{-2} w_k^2 + \gamma_k \log s_k + \beta_k \gamma_k^{-1} \delta_k^{(\mu)}]$$

$$+ e^{-T^{**}u_n} \sum_{k=1}^{n} e^{T^{**}u_k} \beta_k o(\sqrt{\beta_k}) \Bigg).$$

Now, since the sequences $(w_n)$ and $(\delta_n^{(\mu)})$ satisfy Condition (C), the sequences $(\beta_n^2\gamma_n^{-2}w_n^2)$, $(\gamma_n\log s_n)$ and $(\beta_n\gamma_n^{-1}\delta_n^{(\mu)})$ satisfy Condition (C). Moreover, the sequence $(\beta_n^{1/2})$ satisfies Condition (C). The application of Lemma 9 then ensures that, for all $T' \in\,]0, T^{**}[$,

$$\|\Delta_n^{(\theta)}\|_T = O(e^{-T^{**}u_n} + e^{-T'u_n}\mathbb{1}_{b=1} + \beta_n^2\gamma_n^{-2}w_n^2 + \gamma_n\log s_n + \beta_n\gamma_n^{-1}\delta_n^{(\mu)})$$
$$+ o(\sqrt{\beta_n}).$$

Let us recall that $T^{**}$ has been set such that $T^{**} > \mathbb{1}_{b=1}/(2\beta_0)$, and note that $T'$ can be chosen such that $e^{-T^{**}u_n} + e^{-T'u_n} = o(\sqrt{\beta_n})$. Since $\gamma_n\log s_n = o(\sqrt{\beta_n})$, it follows that

$$\|\Delta_n^{(\theta)}\|_T = O(\beta_n^2\gamma_n^{-2}w_n^2 + \beta_n\gamma_n^{-1}\delta_n^{(\mu)}) + o(\sqrt{\beta_n}),$$

which proves (29).

A.4.3. *Proof of* (30). In view of (27), we have

$$\|\Delta_n^{(\mu)}\|_M = O\bigg(e^{-M^*s_n} + e^{-M^*s_n}\sum_{k=1}^n e^{M^*s_k}\gamma_k[\gamma_k\|L_k^{(\mu)}\| + \gamma_k\|R_k^{(\mu)}\| + \|r_k^{(\mu)}\|$$
$$+ \|L_k^{(\theta)}\|^2 + \|R_k^{(\theta)}\|^2 + \|L_k^{(\mu)}\|^2$$
$$+ \|R_k^{(\mu)}\|^2 + \|\Delta_k^{(\theta)}\|_T]\bigg).$$

Replacing $\|r_k^{(\mu)}\|$ by its upper bound given in Assumption (A4)(iv), $\|L_k^{(\theta)}\|$ and $\|L_k^{(\mu)}\|$ by their upper bounds given by Lemma 2, $\|R_k^{(\theta)}\|$ and $\|R_k^{(\mu)}\|$ by their upper bounds given by Lemma 5, $\|\Delta_k^{(\theta)}\|_T$ by its upper bound given in (29), and doing some straightforward simplifications, we deduce that

$$\|\Delta_n^{(\mu)}\|_M = O\bigg(e^{-M^*s_n} + e^{-M^*s_n}\sum_{k=1}^n e^{M^*s_k}\gamma_k[\beta_k^2\gamma_k^{-2}w_k^2 + \gamma_k\log s_k$$
$$+ \beta_k\gamma_k^{-1}\delta_k^{(\mu)}]$$
$$+ e^{-M^*s_n}\sum_{k=1}^n e^{M^*s_k}\gamma_k o(\sqrt{\beta_k})\bigg).$$

Since the sequences $(w_n)$ and $(\delta_n^{(\mu)})$ satisfy Condition (C'), the sequences $(\beta_n^2\gamma_n^{-2}w_n^2)$, $(\gamma_n\log s_n)$, $(n^{-\tilde{s}})$ and $(\beta_n\gamma_n^{-1}\delta_n^{(\mu)})$ satisfy Condition (C'). Moreover, the sequence $(\beta_n^{1/2})$ satisfies Condition (C'). The application of Lemma 10 then ensures that

$$\|\Delta_n^{(\mu)}\|_M = O(e^{-M^*s_n} + \beta_n^2\gamma_n^{-2}w_n^2 + \gamma_n\log s_n + \beta_n\gamma_n^{-1}\delta_n^{(\mu)}) + o(\sqrt{\beta_n}).$$

Noting that $e^{-M^* s_n} + \gamma_n \log s_n = o(\sqrt{\beta_n})$, it follows that

$$\|\Delta_n^{(\mu)}\|_M = O(\beta_n^2 \gamma_n^{-2} w_n^2 + \beta_n \gamma_n^{-1} \delta_n^{(\mu)}) + o(\sqrt{\beta_n}),$$

which concludes the proof of (30).

## REFERENCES

[1] Baras, J. S. and Borkar, V. S. (2000). A learning algorithm for Markov decision processes with adaptive state aggregation. In *Proc. 39th IEEE Conference on Decision and Control.* IEEE, New York.

[2] Benveniste, A., Métivier, M. and Priouret, P. (1990). *Adaptive Algorithms and Stochastic Appproximations.* Springer, Berlin. MR1082341

[3] Bhatnagar, S., Fu, M. C., Marcus, S. I. and Bathnagar, S. (2001). Two timescale algorithms for simulation optimization of hidden Markov models. *IIE Transactions* **3** 245–258.

[4] Bhatnagar, S., Fu, M. C., Marcus, S. I. and Fard, P. J. (2001). Optimal structured feedback policies for ABR flow control using two timescale SPSA. *IEEE/ACM Transactions on Networking* **9** 479–491.

[5] Borkar, V. S. (1997). Stochastic approximation with two time scales. *Systems Control Lett.* **29** 291–294. MR1432654

[6] Delyon, B. and Juditsky, A. B. (1992). Stochastic optimization with averaging of trajectories. *Stochastics Stochastic Rep.* **39** 107–118. MR1275360

[7] Dippon, J. and Renz, J. (1996). Weighted means of processes in stochastic approximation. *Math. Methods Statist.* **5** 32–60. MR1386824

[8] Dippon, J. and Renz, J. (1997). Weighted means in stochastic approximation of minima. *SIAM J. Control Optim.* **35** 1811–1827. MR1466929

[9] Duflo, M. (1996). *Algorithmes Stochastiques.* Springer, Berlin. MR1612815

[10] Ljung, L., Pflug, G. and Walk, H. (1992). *Stochastic Approximation and Optimization of Random Systems.* Birkhäuser, Boston. MR1162311

[11] Horn, R. A. and Johnson, C. R. (1985). *Matrix Analysis.* Cambridge Univ. Press. MR0832183

[12] Konda, V. R. and Borkar, V. S. (1999). Actor-critic like learning algorithms for Markov decision processes. *SIAM J. Control Optim.* **38** 94–123. MR1740605

[13] Konda, V. R. and Tsitsiklis, J. N. (2003). On actor-critic algorithms. *SIAM J. Control Optim.* **42** 1143–1166. MR2044789

[14] Konda, V. R. and Tsitsiklis, J. N. (2004). Convergence rate of linear two-timescale stochastic approximation. *Ann. Appl. Probab.* **14** 796–819. MR2052903

[15] Kushner, H. J. and Clark, D. S. (1978). *Stochastic Approximation Methods for Constrained and Unconstrained Systems.* Springer, New York. MR0499560

[16] Kushner, H. J. and Yang, J. (1993). Stochastic approximation with averaging of the iterates: Optimal asymptotic rate of convergence for general processes. *SIAM J. Control Optim.* **31** 1045–1062. MR1227546

[17] Kushner, H. J. and Yin, G. G. (1997). *Stochastic Approximation Algorithms and Applications.* Springer, New York. MR1453116

[18] Ljung, L. (1978). Strong convergence of a stochastic approximation algorithm. *Ann. Statist.* **6** 680–696. MR0464516

[19] Mokkadem, A. and Pelletier, M. (2005). The compact law of the iterated logarithm for multivariate stochastic approximation algorithms. *Stochastic Anal. Appl.* **23** 181–203. MR2123951

[20] MOKKADEM, A. and PELLETIER, M. (2004). A companion for the Kiefer–Wolfowitz–Blum stochastic approximation algorithm. To appear.

[21] NEVELS'ON, M. B. and HAS'MINSKII, R. Z. (1976). *Stochastic Approximation and Recursive Estimation*. Amer. Math. Soc., Providence, RI. MR0423714

[22] PELLETIER, M. (1998). On the almost sure asymptotic behaviour of stochastic algorithms. *Stochastic Process. Appl.* **78** 217–244. MR1654569

[23] PELLETIER, M. (2000). Asymptotic almost sure efficiency of averaged stochastic algorithms. *SIAM J. Control Optim.* **39** 49–72. MR1780908

[24] POLYAK, B. T. (1990). New method of stochastic approximation type. *Automat. Remote Control* **51** 937–946. MR1071220

[25] POLYAK, B. T. and JUDITSKY, A. B. (1992). Acceleration of stochastic approximation by averaging. *SIAM J. Control Optim.* **30** 838–855. MR1167814

[26] RUPPERT, D. (1991). Stochastic approximation. In *Handbook of Sequential Analysis* (B. K. Ghosh and P. K. Sen, eds.) 503–529. Dekker, New York. MR1174318

[27] YIN, G. (1991). On extensions of Polyak's averaging approach to stochastic approximation. *Stochastics Stochastic Rep.* **33** 245–264. MR1128497

DÉPARTEMENT DE MATHÉMATIQUES
UNIVERSITÉ DE VERSAILLES-SAINT-QUENTIN
45 AVENUE DES ETATS-UNIS
78035 VERSAILLES CEDEX
FRANCE
E-MAIL: mokkadem@math.uvsq.fr
           pelletier@math.uvsq.fr